\documentclass[11pt,reqno]{amsart}
\usepackage{lmodern}
\usepackage[T1]{fontenc}
\usepackage[utf8]{inputenc}
\usepackage{mathrsfs}
\usepackage{mathtools}
\usepackage{enumitem}
\usepackage{amstext}
\usepackage{amsthm}
\usepackage{amssymb}
\usepackage[a4paper]{geometry}
\usepackage[bookmarks=true,bookmarksnumbered=false,bookmarksopen=false,
 breaklinks=false,pdfborder={0 0 0},pdfborderstyle={},backref=false,colorlinks=false]
 {hyperref}
\hypersetup{
 hidelinks}

\makeatletter
\numberwithin{equation}{section}
\usepackage{amsthm}\usepackage{mathrsfs}\usepackage[expansion=false]{microtype}\usepackage{xcolor}
\usepackage[capitalise,nameinlink]{cleveref}
\AtBeginDocument{%
  \@ifundefined{firstaid@cref@updatelabeldata}{}{%
    \gdef\cref@fix@ctr{}\gdef\cref@fix@type{}%
    \let\cref@fix@orig@optarg\refstepcounter@optarg
    \def\refstepcounter@optarg[#1]#2{%
      \gdef\cref@fix@ctr{#2}\gdef\cref@fix@type{#1}%
      \cref@fix@orig@optarg[#1]{#2}}%
    \let\cref@fix@orig@noarg\refstepcounter@noarg
    \def\refstepcounter@noarg#1{%
      \gdef\cref@fix@ctr{#1}\gdef\cref@fix@type{#1}%
      \cref@fix@orig@noarg{#1}}%
    \renewcommand\firstaid@cref@updatelabeldata[1]{%
      \cref@constructprefix{#1}{\cref@result}%
      \edef\cref@fix@tmp{#1}%
      \ifx\cref@fix@tmp\cref@fix@ctr
        \let\@tempa\cref@fix@type
      \else
        \@ifundefined{cref@#1@alias}%
          {\def\@tempa{#1}}%
          {\edef\@tempa{\csname cref@#1@alias\endcsname}}%
      \fi
      \protected@xdef\cref@gcurrentlabel@temp{%
        [\@tempa][\arabic{#1}][\cref@result]%
        \csname p@#1\expandafter\endcsname\csname the#1\endcsname}%
      \aftergroup\firstaid@cref@smugglelabel}%
  }%
}

\crefformat{enumi}{(#2#1#3)}
\Crefformat{enumi}{(#2#1#3)}
\crefrangeformat{enumi}{(#3#1#4)--(#5#2#6)}
\crefmultiformat{enumi}{(#2#1#3)}{ and~(#2#1#3)}{, (#2#1#3)}{ and~(#2#1#3)}
\Crefmultiformat{enumi}{(#2#1#3)}{ and~(#2#1#3)}{, (#2#1#3)}{ and~(#2#1#3)}

\setlist[enumerate]{itemsep=2pt,topsep=4pt}

\theoremstyle{plain}
\newtheorem{theorem}{Theorem}[section]
\newtheorem{proposition}[theorem]{Proposition}\newtheorem{lemma}[theorem]{Lemma}\newtheorem{corollary}[theorem]{Corollary}\theoremstyle{definition}
\newtheorem{definition}[theorem]{Definition}\theoremstyle{remark}
\newtheorem*{remark*}{Remark}\newtheorem{remark}[theorem]{Remark}

\newcommand{\N}{\mathbb{N}}      % the natural numbers
\newcommand{\Z}{\mathbb{Z}}      % the integers
\newcommand{\Q}{\mathbb{Q}}      % the rationals
\newcommand{\R}{\mathbb{R}}      % the reals
\newcommand{\C}{\mathbb{C}}      % the complex numbers
\newcommand{\A}{\mathbb{A}}      % affine space
\newcommand{\Adeles}{\mathbb{A}} % ring of adeles
\newcommand{\Orbit}{\mathscr O}  % conjugacy orbit; O_y denotes its integral
\newcommand{\Phiind}{\Phi_{\mathrm{ind}}} % non-divisible relative roots
\newcommand{\OO}{\mathcal{O}}    % ring of integers of F
\newcommand{\bG}{\mathbf{G}}     % the connected reductive group
\newcommand{\Hk}{\mathcal{H}}    % Hecke algebras
\newcommand{\cO}{\mathcal{C}}    % constructible functions
\newcommand{\gDP}{\mathcal{L}_{\mathrm{gDP}}} % the generalized Denef--Pas language
\newcommand{\VF}{\mathrm{VF}}    % the valued field sort
\newcommand{\VG}{\mathrm{VG}}    % the value group sort
\newcommand{\RF}[1]{\mathrm{RF}_{#1}} % the residue ring sorts
\newcommand{\abs}[1]{\lvert#1\rvert}
\newcommand{\norm}[1]{\lVert#1\rVert}

\DeclareMathOperator{\Span}{span}
\DeclareMathOperator{\vol}{vol}
\DeclareMathOperator{\supp}{supp}

\DeclareMathOperator{\Ad}{Ad}
\DeclareMathOperator{\Hom}{Hom}
\DeclareMathOperator{\val}{val}
\DeclareMathOperator{\im}{im}
\DeclareMathOperator{\Lie}{Lie}
\DeclareMathOperator{\GL}{GL}
\DeclareMathOperator{\Mat}{M}

\DeclareMathOperator{\res}{res}
\DeclareMathOperator{\ac}{ac}

\newcommand{\shelltest}{f} % Test function attached to a shell.
\newcommand{\heckefn}[2]{h_{#1,#2}} % Hecke function indexed by level and group element.
\newcommand{\levsupport}{Levi support hypothesis}
\newcommand{\matcong}[2]{K_{#1}^{#2}}
\newcommand{\matball}[1]{B_{#1}}
\newcommand{\ball}{B}
\newcommand{\adball}[1]{B_{\mathrm{Ad},#1}}
\newcommand{\adnorm}[2]{\norm{#2}_{\mathrm{Ad},#1}}
\newcommand{\Hsc}{\Hk_{\mathrm{sc}}}
\newcommand{\bS}{\mathbf S}
\newcommand{\bH}{\mathbf H}
\newcommand{\adtor}{\mathbf A}
\newcommand{\centrallattice}{\Lambda_Z}
\newcommand{\Gbar}{\overline{G}}
\newcommand{\Irrsc}[1]{\operatorname{Irr}_{\mathrm{sc}}(#1)}
\newcommand{\IrrscZ}[1]{\Irrsc{#1;\centrallattice}}
\newcommand{\Loc}{\mathrm{Loc}_0}
\DeclareMathOperator{\SL}{SL}
\newtheorem{example}[theorem]{Example}\crefformat{enumii}{(#2#1#3)}
\Crefformat{enumii}{(#2#1#3)}

\makeatother

\begin{document}

\title[Support bounds and Bernstein stabilization]{Linear support bounds for
supercuspidal matrix coefficients of $p$-adic groups}
\author[I. Glazer]{Itay Glazer}
\author[E. Lapid]{Erez Lapid}
\author[Z. Mao]{Zhengyu Mao}
\author[M. Shusterman]{Mark Shusterman}
\date{November 18, 2026}
\subjclass[2020]{22E50, 03C98, 11F70}
\keywords{Supercuspidal representations, matrix coefficients, Bernstein stabilization,
second adjointness, constructible functions, limit multiplicities}
\begin{abstract}
Let $\bG$ be a connected reductive group over a non-archimedean local
field $F$ of characteristic zero. For any standard congruence tower
$K_{\bullet}=(K_{n})$ in $\bG(F)$, supercuspidal matrix coefficients
between $K_{n}$-fixed vectors are supported in an adjoint ball of
radius $an+b_{K_{\bullet}}$, compact modulo the center, with $a$
depending only on $F$ and $\bG$. When $\bG(F)$ is non-compact modulo
its center, the linear dependence on $n$ is optimal. In algebraic
families we obtain a uniform radius $An+Bk+C+\val_{F}(N)$, where
$k$ bounds the logarithmic height of the parameter.
We deduce a sharp linear bound for the exponent in Bernstein's stabilization theorem.
For quasi-split classical groups, we also remove the dyadic level
restriction in the limit multiplicity theorem.
\end{abstract}

\maketitle
\tableofcontents{}

\section{Introduction}

\label{sec:intro} Unless otherwise specified, $F$ is a non-archimedean
local field of characteristic zero, with normalized absolute value
$|\cdot|_{F}$, ring of integers $\OO_{F}$, maximal ideal $\mathfrak{p}_{F}$,
residue field of cardinality $q_{F}$, and $\operatorname{val}_{F}$
for its integer-normalized valuation. When $F$ is fixed, we may omit
the subscript $F$. Let $\bG$ be a connected reductive group over
$F$, let $G=\bG(F)$, and let $Z$ be its center. Representations
are smooth and complex. We write $V^{\vee}$ for the smooth dual of
$V$. We denote by $\Irrsc G$ the set of isomorphism classes of irreducible
supercuspidal representations of $G$.

Fix a standard tower $K_{\bullet}:=(K_{n})_{n\ge1}$ in $G$ as in
\cref{subsec:towers-length}. The matrix congruence subgroups $\matcong n{\rho}$
attached to a closed embedding $\rho:\bG\hookrightarrow\GL_{d}$ give
an example after a fixed shift of level. For $R\ge0$, let 
\[
\ball(R)=\{g\in G:gK_{j+\lfloor R\rfloor}g^{-1}\subseteq K_{j}\text{ for all sufficiently large }j\},
\]
the adjoint ball of radius $R$ for this tower. These balls are invariant
under multiplication by $Z$ and compact modulo $Z$. Changing the
standard tower changes the radii by at most an additive constant (\cref{subsec:towers-length}).
Our main result is the following.
\begin{theorem}[Linear support bound]
\label{thm:main} There exists $a>0$, depending only on $F$ and
$\bG$, such that for every standard tower $K_{\bullet}=(K_{n})$,
there exists $b_{K_{\bullet}}\ge0$ with the following property. For
every $\pi\in\Irrsc G$, every $n\ge1$, and every $v\in\pi^{K_{n}}$,
$v^{\vee}\in(\pi^{\vee})^{K_{n}}$, 
\[
\supp\bigl(g\mapsto\langle\pi(g)v,v^{\vee}\rangle\bigr)\subset\ball(an+b_{K_{\bullet}}).
\]
\end{theorem}

\cref{thm:main} amounts to property (PSC), introduced by Finis--Lapid--Müller
\cite{MR3001800}, for every connected reductive group over $F$ (\cref{lem:psc}).
There is no restriction on the residual characteristic, the ramification
of $F$, or the splitting of $\bG$.

\cref{thm:main} is trivial when $G/Z$ is compact, in which case
$G$ is a single adjoint ball. Otherwise the linear order in \cref{thm:main}
is optimal. More precisely, for every $\pi\in\Irrsc G$, the smallest
$R$ for which $\ball(R)$ contains the supports of all matrix coefficients
of $\pi$ between $K_{n}$-fixed vectors is at least $2n-O_{\pi,K_{\bullet}}(1)$
as $n\to\infty$ (\cref{rem: sharpbnd}). In particular, the slope
$a$ in \cref{thm:main} satisfies $a\ge2$.

For each fixed level, Bernstein's finiteness theorem gives a support
bound uniform in the supercuspidal representation; see \cite[Remark~14]{MR3001800}.
\cref{thm:main} shows that this bound can be chosen to grow linearly
with the level. For $\GL_{d}$ over every non-archimedean local field,
in either characteristic, the bound for the usual matrix norm and
principal congruence subgroups was proved in \cite{MR3930019} using
Rankin--Selberg integrals, with an explicit constant depending only
on $d$.

A different route to (PSC) is available when every supercuspidal representation
is contained in a representation induced from a cuspidal representation
of an open subgroup compact modulo the center \cite[Corollary~13]{MR3001800}.
The constants in this support bound are explicit \cite[Lemma~8 and Corollary~10]{MR3001800}.
See also \cite[Theorem~3]{MR3940435}. Such classification results
are known in many cases. They include $\GL_{d}$ and $\SL_{d}$ through
Bushnell--Kutzko's work \cite{MR1204652,MR1209709}, classical groups
in odd residual characteristic through Stevens' work \cite{MR2390287},
and the almost-all-places statement for a group over a number field
through Kim's exhaustion theorem \cite{MR2276772} for Yu's construction
of tame supercuspidal representations \cite{MR1824988}. Fintzen \cite{MR4199732}
proved exhaustion for tamely ramified groups whenever the residue
characteristic does not divide the order of the absolute Weyl group.
Our proof uses no classification or construction of supercuspidal
representations.

The adjoint support bound of \cref{thm:main} is uniform in algebraic
families in the following sense. Let $\bS$ be a closed subscheme
of $\A^{D}_{\Z}$ and let $\mathscr{G}$ be a reductive group scheme
over $\bS$, i.e., a smooth affine group scheme with connected reductive
geometric fibers. Fix a closed immersion $\rho:\mathscr{G}\hookrightarrow\GL_{d,\bS}$.
For $s\in\bS(F)$, put $G_{s}=\mathscr{G}_{s}(F)$ and let $\ball_{s}(R)$
be the adjoint ball in $G_{s}$ attached to the tower $(\matcong n{\rho_{s}})_{n\ge1}$
as above.
\begin{theorem}[Support in algebraic families]
\label{thm:families-intro} There exist integers $A,B,C\ge0$ and
$N\ge1$, depending only on $(\mathscr{G},\bS,\rho)$, such that the
following holds over every non-archimedean local field $F$ of characteristic
zero. For every $k\in\Z_{\ge0}$, every $s\in\bS(F)$ with $\max_{1\le i\le D}|s_{i}|_{F}\le q^{k}_{F}$,
and every $n\ge1$, the matrix coefficients between $\matcong n{\rho_{s}}$-fixed
vectors of every irreducible supercuspidal representation of $G_{s}$
are supported in 
\[
\ball_{s}\bigl(An+Bk+C+\val_{F}(N)\bigr).
\]
\end{theorem}

In fact, the proof gives the same assertion with $Z_{s}\matball{\rho_{s}}(R)$
in place of $\ball_{s}(R)$. Here $Z_{s}=Z(G_{s})$ and $\matball{\rho_{s}}(R)=\{g\in G_{s}:\norm g_{\rho_{s}}\le q^{R}_{F}\}$,
where $\norm g_{\rho_{s}}$ is the largest absolute value of the entries
of $\rho_{s}(g)$ and its inverse. For the Lie lattice induced by
$\rho_{s}$, \eqref{eq:matrix-adjoint-balls} gives $Z_{s}\matball{\rho_{s}}(R)\subseteq\ball_{s}(2R)$.
Thus, if $A,B,C,N$ are the constants in the matrix bound, then $2A,2B,2C,N^{2}$
work in \cref{thm:families-intro}.

For integral parameters, i.e., those for which $\norm s_{F}:=\max_{1\le i\le D}|s_{i}|_{F}\le1$,
one may take $k=0$. The proof establishes the existence of these
constants without giving numerical values for them.

Fix a connected reductive group $\bG$ over a number field $k$ and a
closed embedding $\rho:\bG\hookrightarrow\GL_{d,k}$.
\begin{corollary}[Uniformity over a number field, \cref{cor:rational-model}]
\label{thm:rational-intro} There exist $A,B\ge0$ and an integer $N\ge1$,
depending only on $(k,\bG,\rho)$, such that for every finite place $v$ of
$k$, every finite extension $F/k_{v}$, and every $n\ge1$, the matrix
coefficients between $\matcong n{\rho}$-fixed vectors in every irreducible
supercuspidal representation of $\bG(F)$ are supported in
\[
Z(\bG(F))\matball{\rho}\bigl(An+B+\val_{F}(N)\bigr).
\]
\end{corollary}

The valuation term vanishes outside the finitely many residue characteristics
dividing $N$. At the remaining primes it accounts for ramification. 

For a standard tower $K_{\bullet}=(K_{n})$, the support bound of
\cref{thm:main} gives an upper bound for the size of the supercuspidal
part $\Hsc(G,K_{n})$ of the Hecke algebra with respect to $K_{n}$.
For simplicity we state it here when the center of $G$ is compact.
Upon choosing a Haar measure on $G$, we identify $\Hsc(G,K_{n})$
with the span of the $K_{n}$-biinvariant matrix coefficients of irreducible
supercuspidal representations of $G$.
\begin{corollary}[The supercuspidal part of the Hecke algebra; see \cref{cor:hecke-upper}]
\label{cor:hecke-intro} Suppose that the center of $G$ is compact
and let $r$ be the $F$-rank of $G$. There is a constant $C>0$,
depending only on $F$, $G$, and the chosen tower, such that 
\[
\dim\Hsc(G,K_{n})=\sum_{\pi\in\Irrsc G}(\dim\pi^{K_{n}})^{2}\le C(1+n)^{r}q^{n\dim G}\qquad(n\ge1).
\]
\end{corollary}

Under the same assumption, \cite[Theorem~6.14]{LS26} gives, for all
sufficiently large $n$, 
\[
q^{n\dim G}\ll\dim\Hsc(G,K_{n}).
\]
Thus the upper and lower bounds differ by a factor of order of magnitude $(1+n)^{r}$.

Property (BD) (bounded degree), introduced by Finis--Lapid--Müller
\cite{MR3352530}, concerns standard unnormalized intertwining operators
for maximal parabolic subgroups of Levi subgroups of $G$. The matrix
entries of these operators on spaces of congruence-fixed vectors are
rational functions of $q^{-s}$, where $s$ is the one-dimensional
unramified twisting parameter. The denominator degrees are bounded
uniformly in terms of $G$ only. In particular, property (BD) implies that the numerator
degrees on vectors of congruence level $n$ are at most $Cn$, uniformly
in the irreducible inducing representation. The implication from (PSC)
to (BD) of Finis--Lapid \cite{MR3940435} gives the following consequence. 
\begin{corollary}[Local degree bounds]
\label{cor:bd-intro} Every connected reductive group over $F$ satisfies
property (BD). 
\end{corollary}

Over a number field, the argument in \cref{sec:limit} gives (BD)
uniformly at all finite places. Combining this with the limit-multiplicity results
of \cite{MR3940435} yields the following global application. 
\begin{corollary}[Global degree bounds and limit multiplicity]
 \label{cor:global-intro}Let $\bH$ be a connected reductive group
over a number field $k$. Then $\bH$ satisfies property (BD) uniformly
over all finite places. If $\bH$ also satisfies (TWN), then, for
every finite set $S$ of places containing all archimedean places,
the limit multiplicity property holds for every non-degenerate family
of open subgroups of a fixed open compact subgroup of $\bH(\Adeles^{S})$. 
\end{corollary}

The definitions of (BD), (TWN), and non-degeneracy are those of \cite{MR3940435,MR3352530};
see \cref{cor:limit,rem:all-places}. For quasi-split classical groups,
the limit-multiplicity assertion removes the bounded-level assumption
at the dyadic places outside $S$ in \cite[Corollary~2(2)]{MR3940435}.
At places of odd residue characteristic, our proof of (BD) also avoids
the use of Stevens's exhaustion theorem for supercuspidal representations.

We now sketch the proof of \cref{thm:main,thm:families-intro}. The
proof has a representation-theoretic component and a model-theoretic
component.

For the first component, fix a local field $F$, a group $G=\bG(F)$,
and a closed embedding $\rho:\bG\hookrightarrow\GL_{d}$. Suppose
for simplicity that $Z$ is compact and $\rho(Z)\subseteq\GL_{d}(\OO_{F})$.
Fix a Haar measure $dg$ on $G$ and identify Hecke measures with
their densities. For each $g\in G$, define the $\matcong n{\rho}$-biinvariant
function $\heckefn ng=1_{\matcong n{\rho}}*\delta_{g}*1_{\matcong n{\rho}}$
on $G$, where $1_{\matcong n{\rho}}$ is the unnormalized indicator.
The Hecke measure $\shelltest_{n,r}$ is given by the shell integral
\[
\shelltest_{n,r}=\int_{\norm g_{\rho}=q^{r}_{F}}\heckefn ng*\heckefn n{g^{-1}}\,dg\ \in\ \Hk(G,\matcong n{\rho}).
\]
For any admissible representation $\pi$ with a unitary realization,
the trace $\Theta_{\pi}(\shelltest_{n,r})$ is the integral over the
shell $\norm g_{\rho}=q^{r}_{F}$ of the squared Hilbert--Schmidt
norm of $\pi(\heckefn ng)$ (see \eqref{eq:HS}). Since $\pi(1_{\matcong n{\rho}})$
is a positive multiple of the orthogonal projection onto $\pi^{\matcong n{\rho}}$,
this integral vanishes precisely when all matrix coefficients between
$\matcong n{\rho}$-fixed vectors vanish on that shell. If $\pi$
is irreducible supercuspidal, then the integral vanishes for all sufficiently
large $r$. The problem is to bound this threshold uniformly in $\pi$.

For fixed $n$, Clozel's proof of Howe's finiteness conjecture \cite{MR986793}
shows that the sequences $r\mapsto O_{y}(\shelltest_{n,r})$, with
$y$ ranging over the compact regular semisimple elements of $G$,
span a finite-dimensional space. Weyl's integration formula and Deligne's
support theorem for supercuspidal characters \cite{MR0425033} place
$r\mapsto\Theta_{\pi}(\shelltest_{n,r})$ in this space, as shown
in \cref{lem:fixed-level-characters}. The eventually vanishing sequences in this span
form a finite-dimensional subspace. The union of the supports of the
elements of a finite basis is finite. Thus, for fixed $G,\rho$, there
is a threshold $R(n)$ such that $\Theta_{\pi}(\shelltest_{n,r})=0$
for every irreducible supercuspidal $\pi$ and every $r\ge R(n)$.

The model-theoretic component controls this threshold as the field
and the group vary, and is based on theory of motivic integration
developed by Cluckers--Loeser \cite{MR2403394,MR2630060}, extended to all
characteristic zero local fields by Cluckers--Halupczok \cite{MR3732692}
and Cluckers--Gordon--Halupczok \cite{MR3859937}. 

We realize $(\bG,\rho)$ as a fiber of an embedded algebraic family
$\mathscr{G}\hookrightarrow\GL_{d,\bS}$ and write $\shelltest_{s,n,r}$
for the test function on the fiber $G_{s}$, $s\in\bS(F)$. With the
Haar and orbital measures chosen as in \cref{sec:orbital}, the orbital
integrals form a constructible family in $(s,y,n,r)$, uniformly in
$F$, as in the framework of \cite{MR3732692,MR3859937}. 

The proof of \cref{prop:threshold} gives a definable threshold $R_{F}(s,n)$
such that, for $r\ge R_{F}(s,n)$, each sequence $r\mapsto O_{y}(\shelltest_{s,n,r})$
agrees with an element of a finite-dimensional space $E_{F}$ of exponential
polynomials on congruence classes, with congruence modulus, exponents,
and degree bounds independent of $F,s,n$ and $y$. Moreover, 
\[
R_{F}(s,n)\le An+\gamma_{F}(s),
\]
where $\gamma$ is a definable integer-valued function. In particular,
the character sequence $r\mapsto\Theta_{\pi}(\shelltest_{s,n,r})$
agrees with an element of $E_{F}$ for $r\geq R_{F}(s,n)$. Since
$\Theta_{\pi}(\shelltest_{s,n,r})$ vanishes for $r$ large and an
exponential polynomial which vanishes eventually vanishes identically,
it vanishes for all $r\ge An+\gamma_{F}(s)\geq R_{F}(s,n)$. This
gives the support bound for the matrix tower. On a fixed fiber, passage
to any standard tower changes only the additive constant in the adjoint
bound, proving \cref{thm:main}.

The uniform bound in \cref{thm:families-intro} comes from an additional
observation. Conjugation by $\GL_{d}(\OO_{F})$ preserves both the
matrix congruence subgroups and the matrix norm, so a support bound
for one embedded fiber applies to every integrally conjugate fiber.
By local rigidity of embeddings of reductive groups (\cref{lem:local-integral-conjugacy}),
for every $s'$ sufficiently close to $s$, the embedded fiber $\rho_{s'}(G_{s'})$
is $\GL_{d}(\OO_{F})$-conjugate to the embedded fiber $\rho_{s}(G_{s})$.
Thus, for fixed $n$, the minimum of $R_{F}(s',n)$ over all $s'$ for which
$\rho_{s'}(G_{s'})$ is integrally conjugate to $\rho_{s}(G_{s})$ yields
a common threshold that is locally constant in $s$. Taking its maximum on the compact parameter set $\norm s_{F}\le q^{k}_{F}$
gives a definable function of $k$ and $n$. \cref{lem:vg-bound}
then bounds this threshold by $An+Bk+C+\val_{F}(N)$, which proves
\cref{thm:families-intro}.

The second part of the paper deduces quantitative stabilization from
support bounds. In the stabilization statements below and in \crefrange{sec:strata}{sec:adjunction},
$K_{n}$ denotes the subgroup of level $n$ in the adapted tower of
\cref{subsec:adapted-filtrations}. In characteristic zero, \cref{thm:main}
gives a linear adjoint support bound on each Levi subgroup $L$ for
the induced tower $(L\cap K_{n})$, which is standard after a fixed
shift (see \cref{subsec:adapted-filtrations,sec:averaging}). The
deductions from these bounds hold in any characteristic.

Fix a proper standard parabolic $P$ with standard Levi decomposition
$M\ltimes U$. Let $A_{M}$ be the group of $F$-points of the maximal
$F$-split torus $\mathbf{A}_{M}$ in the center of $M$. For $a\in A_{M}$
define its contraction rate by 
\[
\epsilon_{P}(a)=\min_{\alpha\in\Phi(A_{M},U)}\val_{F}(\alpha(a)),
\]
where $\Phi(A_{M},U)$ is the set of weights of $A_{M}$ in $\Lie U$.
Thus $\epsilon_{P}(a)>0$ means that conjugation by $a$ contracts
$U$. Let $e_{K_{n}}$ denote normalized averaging over $K_{n}$,
so that $\pi(e_{K_{n}})$ is a projection of $V$ on $V^{K_{n}}$,
and set 
\[
T_{n,a}=\pi(e_{K_{n}}*\delta_{a}*e_{K_{n}})=\pi(e_{K_{n}})\pi(a)\pi(e_{K_{n}})\text{ on }V^{K_{n}}.
\]

\begin{theorem}[Stabilization for sufficiently contracting elements; \cref{thm:stab}]
 \label{thm:stab-intro} There exists $c>0$ such that for every
smooth (not necessarily admissible) representation $V$, every $n\ge1$,
and every $a\in A_{M}$ with $\epsilon_{P}(a)\ge cn$, 
\[
V^{K_{n}}=\ker T_{n,a}\oplus\im T_{n,a}.
\]
For fixed $P$ and $n$, the two summands are independent of $a$
in this range and $T_{n,a}$ is invertible on its image. The Jacquet
map $j_{U}:V\to V_{U}$ vanishes on $\ker T_{n,a}$ and restricts
to an isomorphism 
\[
\im T_{n,a}\xrightarrow{\;\sim\;}(V_{U})^{K_{n}\cap M}.
\]
\end{theorem}

In particular, for a contracting $t\in A_{M}$ we have $T_{n,t^{r}}=T^{\,r}_{n,t}$
and $\epsilon_{P}(t^{r})=r\epsilon_{P}(t)$ for any integer $r\ge1$.
Hence, the powers of $T_{n,t}$ stabilize once $r\epsilon_{P}(t)\ge cn$
(\cref{cor:stab-powers}). Qualitative stabilization is due to Bernstein
\cite{BernsteinHarvard}. Our assertion is the linear dependence on
the congruence level. The support bound first controls averaging in
a supercuspidal inducing representation. Applying this successively
along the Bruhat filtration bounds the annihilation of the Jacquet
kernel in an induced representation; an algebraic reduction then passes
the bound to all smooth representations. The precise formulation and
proof are given in \cref{sec:stab}.

Let $\bar{P}=M\bar{U}$ be opposite to $P$, and put $\Delta_{P}(m)=|\det(\Ad(m)|_{\Lie U})|$
and $r_{P}(V)=\Delta^{-1/2}_{P}\otimes V_{U}$. Write $j_{U}$ and
$j_{\bar{U}}$ for the Jacquet maps.

Suppose that $(\pi,V)$ is admissible. Then for every $v\in V$ and
$v^{\vee}\in V^{\vee}$ the germ of the normalized matrix coefficient
$\Delta_{P}(a)^{-1/2}\langle\pi(a)v,v^{\vee}\rangle$ as $a\to0$
in the contracting cone is $A_{M}$-finite and depends only on $j_{U}(v)$
and $j_{\bar{U}}(v^{\vee})$. More precisely, it is represented by
the $A_{M}$-finite function 
\[
B_{P}\bigl(r_{P}(\pi)(a)j_{U}(v),j_{\bar{U}}(v^{\vee})\bigr),\qquad a\in A_{M},
\]
where $B_{P}$ is the Casselman pairing between $r_{P}(V)$ and $r_{\bar{P}}(V^{\vee})$.
Recall that smooth $A_{M}$-finite functions are finite linear combinations
of functions $a\mapsto\chi(a)p(\nu_{A_{M}}(a))$, where $\chi$ is
a smooth character of $A_{M}$ and $p$ is a complex polynomial on
the cocharacter lattice $X_{*}(\mathbf{A}_{M})$ and $\nu_{A_{M}}:A_{M}\to X_{*}(\mathbf{A}_{M})$
is the valuation homomorphism. Even if $\pi$ is not admissible (but
only smooth), the relation 
\begin{equation}
\Delta_{P}(a)^{-1/2}\langle\pi(a)v,v^{\vee}\rangle=B_{P}\bigl(r_{P}(\pi)(a)j_{U}(v),j_{\bar{U}}(v^{\vee})\bigr)\label{eq: asymp}
\end{equation}
still holds for $a\rightarrow0$ and defines the Casselman pairing.

%Here $A_M$-finite means that the translates span a finite-dimensional space.
\begin{corollary}[A linear threshold for matrix-coefficient asymptotics]
 \label{cor:pairing-intro} For every smooth representation $(\pi,V)$,
the relation \eqref{eq: asymp} holds for every $n\ge1$, $v\in V^{K_{n}}$,
$v^{\vee}\in(V^{\vee})^{K_{n}}$, and $a\in A_{M}$ with $\epsilon_{P}(a)\ge cn$.
Here $c$ is the constant of \cref{thm:stab-intro}. 
\end{corollary}

 The qualitative statements are classical \cite{CassNotes,BernsteinHarvard,MR3430373,MR2567785,ChanDual}.
The assertion here is the linear threshold in the congruence level;
see \cref{lem:lift,thm:adj}. The same threshold enters the finite-level
realization of second adjointness in \cref{sec:adjunction}.

The paper is organized as follows. \Cref{sec:prelim} introduces standard
towers, adjoint balls, and matrix balls. \Cref{sec:constructible,sec:orbital,sec:char0}
prove the support theorem. \Cref{sec:limit} introduces the adapted
tower and gives the applications to the Hecke algebra, degree bounds,
and limit multiplicities. \Cref{sec:strata} treats the Bruhat strata.
\Cref{sec:averaging,sec:smooth,sec:stab} deduce averaging and stabilization
from support bounds, and \cref{sec:pairing,sec:adjunction} give the
applications to canonical lifts, asymptotics, and second adjointness.

\subsection*{AI disclosure}
We used AI tools, including Claude (Anthropic) and Codex (OpenAI), to help
develop and check arguments, locate references, and draft and revise
the manuscript. The authors take full responsibility for its contents.

\subsection*{Acknowledgments}
The work was supported in part by a grant of access to OpenAI models
through the ChatGPT for Academic Researchers Program.
Itay Glazer was supported by ISF grant 3422/24.
Mark Shusterman is supported by the European Union (ERC, Function Fields,
101161909) and is the Dr.~A.~Edward Friedmann Career Development Chair
in Mathematics.

\section{Groups, filtrations, and lengths}

\label{sec:prelim}

For any totally disconnected locally compact group $H$, let $\Hk(H)$
be the convolution algebra of compactly supported locally constant
complex measures on $H$. A choice of left Haar measure $dh$ identifies
$C^{\infty}_{c}(H)$ with $\Hk(H)$ by $f\mapsto f\,dh$. For a compact
open subgroup $J\subset H$, let $e_{J}$ be Haar probability on $J$,
viewed as a measure on $H$. It satisfies $e_{J}*e_{J}=e_{J}$ and
has density $\mathbf{1}_{J}/\vol(J)$ with respect to $dh$. The algebra
of bi-$J$-invariant Hecke measures is 
\[
\Hk(H,J)=e_{J}\Hk(H)e_{J}.
\]
A point mass $\delta_{h}$ acts on $\Hk(H)$ by convolution. For the
reductive groups and their Levi subgroups considered below, left Haar
measure is also right invariant.

\subsection{Standard towers, adjoint balls, and matrix balls}

\label{subsec:towers-length} Choose a uniformizer $\varpi$ of $F$
and normalize $\val$ by $\val(\varpi)=1$. Put $m=\dim\bG$ and $\mathfrak{g}=\Lie(\bG)(F)$.
Following Serre \cite[Part~II, Ch.~IV, \S\S8--9]{MR2179691}, choose
a standard open subgroup $K_{1}\subset G$ and analytic coordinates
$\iota\colon\mathfrak{p}\OO^{m}\to K_{1}$ in which multiplication
is given by a formal group law $\Phi$ over $\OO$. The congruence
subgroups 
\[
K_{n}=\iota(\mathfrak{p}^{n}\OO^{m}),\qquad n\ge1,
\]
form a \emph{standard tower}. Its Lie lattice is $\mathcal{L}=d\iota_{0}(\OO^{m})\subset\mathfrak{g}$.

For an analytic map $u$ fixing zero, Taylor expansion gives $\varpi^{-j}u(\varpi^{j}X)\to du_{0}(X)$
uniformly on any fixed lattice $L$. Thus $u(\varpi^{j}L)\subseteq\varpi^{j}L'$
for all large $j$ if and only if $du_{0}(L)\subseteq L'$, since
the lattice $L'$ in the target is open and closed. Applied to a change
of coordinates and its inverse, this shows that two standard towers
with the same Lie lattice agree at all sufficiently large levels.

In characteristic zero, $K_{n}=\exp(\mathfrak{p}^{n}\mathcal{L})$ for all
sufficiently large $n$. The construction of the tower by a formal group law
also applies in positive characteristic.

Write $\norm{\cdot}_{\mathcal{L}}$ for the operator norm on $\operatorname{End}_{F}(\mathfrak{g})$
induced by the maximum norm in any $\OO$-basis of $\mathcal{L}$.
The adjoint ball of radius $R\ge0$ depends only on $\mathcal{L}$
and is given by 
\[
\ball(R)=\{g\in G:\norm{\Ad(g)}_{\mathcal{L}}\le q^{R}\}.
\]
These balls are invariant under multiplication by $Z$ and compact
modulo $Z$ \cite[\S2]{MR3940435}. For $r=\lfloor R\rfloor$, the
norm condition says $\varpi^{r}\Ad(g)\mathcal{L}\subseteq\mathcal{L}$.
Applying the same Taylor estimate to conjugation by $g$ in the coordinates
$\iota$ gives 
\begin{equation}
g\in\ball(R)\quad\Longleftrightarrow\quad gK_{j+r}g^{-1}\subseteq K_{j}\quad\text{for all sufficiently large }j.\label{eq:adjoint-tower-balls}
\end{equation}
Let $(K'_{j})$ be another standard tower with balls $\ball'(R)$.
The two Lie lattices are commensurable. Applying the estimate to the
coordinate changes in both directions gives an integer $d\ge0$ such
that $K_{j+d}\subseteq K'_{j}$ and $K'_{j+d}\subseteq K_{j}$ for
all large $j$. Applying these inclusions on both sides of \eqref{eq:adjoint-tower-balls}
gives 
\[
\ball(R)\subseteq\ball'(R+2d),\qquad\ball'(R)\subseteq\ball(R+2d).
\]

A closed immersion of algebraic $F$-groups $\rho:\bG\hookrightarrow\GL_{d}$
defines the lattice $\mathcal{L}_{\rho}=(d\rho)^{-1}(\Mat_{d}(\OO))$
and the matrix tower 
\[
\matcong n{\rho}=\{g\in G:\rho(g)\in1+\mathfrak{p}^{n}\Mat_{d}(\OO)\},\qquad n\ge1.
\]
These are compact open normal subgroups of $\matcong 0{\rho}=\rho^{-1}(\GL_{d}(\OO))$
and form a neighborhood basis of the identity. In local coordinates
$X\in\mathfrak{g}$ centered at $1$ with identity differential, $\rho(g)-1=d\rho(X)+O(|d\rho(X)|^{2})$
for the entrywise maximum norm $|\cdot|$. Hence $|\rho(g)-1|=|d\rho(X)|$
near zero. Thus, for all sufficiently large $n$, the coordinate map identifies
$\matcong n{\rho}$ with $\mathfrak{p}^{n}\mathcal{L}_{\rho}$. Rescaling a basis of $\mathcal{L}_{\rho}$ by
$\varpi^{b}$ for $b$ sufficiently large makes the group law integral
and realizes $(\matcong{n+b}{\rho})_{n\ge1}$ as a standard tower
with Lie lattice $\mathfrak{p}^{b}\mathcal{L}_{\rho}$.

On $\GL_{d}(F)$ put 
\[
\norm h_{\mathrm{mat}}=\max_{i,j}(|h_{ij}|,|(h^{-1})_{ij}|)\ge1.
\]
Its pullback is $\norm g_{\rho}=\norm{\rho(g)}_{\mathrm{mat}}$. Define
\[
\ell_{\rho}(g)=\log_{q}\norm g_{\rho}\ge0,\qquad\matball{\rho}(R)=\{g\in G:\ell_{\rho}(g)\le R\}\quad(R\ge0).
\]
The integer-valued length $\ell_{\rho}$ is symmetric and subadditive,
and bi-$\matcong 0{\rho}$-invariant. Define 
\[
\ell_{\rho,Z}(g)=\min_{z\in Z}\ell_{\rho}(zg).
\]
The sets $\matball{\rho}(R)$ are compact and open, and $\ell_{\rho,Z}$
is bi-$\matcong 0{\rho}$-invariant.

In the norm induced by $\mathcal{L}_{\rho}$, conjugation $X\mapsto\rho(g)X\rho(g)^{-1}$
gives $\norm{\Ad(g)}_{\mathcal{L}_{\rho}}\le\norm g^{2}_{\rho}$.
Choose an integer $e\ge0$ with $\varpi^{e}\mathcal{L}_{\rho}\subseteq\mathcal{L}\subseteq\varpi^{-e}\mathcal{L}_{\rho}$.
The operator norm for $\mathcal{L}$ is at most $q^{2e}$ times that
for $\mathcal{L}_{\rho}$. Thus 
\begin{equation}
Z\matball{\rho}(R)\subseteq\ball(2R+2e)\qquad(R\ge0),\label{eq:matrix-adjoint-balls}
\end{equation}
with $e=0$ for the induced lattice $\mathcal{L}=\mathcal{L}_{\rho}$.
Since the matrix tower is standard after a fixed shift, the comparison
of standard towers above gives an integer $\delta\ge0$ such that
\begin{equation}
\matcong{n+\delta}{\rho}\subset K_{n},\qquad K_{n+\delta}\subset\matcong n{\rho}\qquad(n\ge1),\label{eq:standard-matrix-comparison}
\end{equation}
after increasing $\delta$ to cover the finitely many initial levels.

Let $X^{*}_{F}(\bG)$ be the lattice of $F$-rational algebraic characters
of $\bG$. The valuation homomorphism $\nu_{G}\colon G\to\Hom(X^{*}_{F}(\bG),\Z)$,
given by $\nu_{G}(g)(\chi)=\val(\chi(g))$, has open normal kernel
$G^{1}=\bigcap_{\chi\in X^{*}_{F}(\bG)}\ker|\chi|$. Restriction to
the maximal $F$-split central torus $\mathbf{A}_{G}$ identifies
$X^{*}_{F}(\bG)$ with a finite-index subgroup of $X^{*}(\mathbf{A}_{G})$.
Hence $\nu_{G}(\mathbf{A}_{G}(F))$ has finite index in $\Hom(X^{*}_{F}(\bG),\Z)$,
and so does $\nu_{G}(G)$. In particular $G/G^{1}$ is a lattice
of rank $\dim\mathbf{A}_{G}$. Since $\mathbf{A}_{G}(F)$ is cocompact
in $Z$ and has compact valuation kernel, $Z\cap G^{1}$ is compact.
Every compact subgroup of $G$ is contained in $G^{1}$.
\begin{remark}[Sharpness of the support bound]
\label{rem: sharpbnd} Suppose that $G/Z$ is noncompact. For every
$\pi\in\Irrsc G$ and every standard tower $K_{\bullet}=(K_{n})$,
there is a constant $C_{\pi,K_{\bullet}}$ such that, for all sufficiently
large $n$, some matrix coefficient between $K_{n}$-fixed vectors
is nonzero at an element of adjoint radius at least $2n-C_{\pi,K_{\bullet}}$.

Choose a noncentral $F$-cocharacter $\lambda$ of $\bG$ and put
$t=\lambda(\varpi)$. By the additive comparison of levels and radii
above, it suffices to prove the claim for a tower whose Lie lattice
$\mathcal{L}$ is spanned by weight vectors of $\lambda$, since fixed
shifts change only $C_{\pi,K_{\bullet}}$. The extreme weights of
$\lambda$ on $\mathfrak{g}$ are $\pm D$ for some $D>0$, so $\norm{\Ad(t^{\pm1})}_{\mathcal{L}}=q^{D}$.
Choose $u\in\pi$ and $u^{\vee}\in\pi^{\vee}$ with $\langle u,u^{\vee}\rangle=1$.
By smoothness and \eqref{eq:adjoint-tower-balls}, after a fixed shift
of the tower, $K_{1}$ fixes both vectors and $t^{\pm1}K_{j+D}t^{\mp1}\subseteq K_{j}$
for every $j\ge1$.

For $m=\lfloor(n-1)/D\rfloor$, iteration gives $t^{\pm m}K_{n}t^{\mp m}\subseteq K_{1}$.
Thus $\pi(t^{-m})u$ and $\pi^{\vee}(t^{m})u^{\vee}$ are $K_{n}$-fixed,
and their matrix coefficient at $t^{2m}$ equals $\langle u,u^{\vee}\rangle=1$.
The weight $-D$ gives the eigenvalue $\varpi^{-2mD}$ of $\Ad(t^{2m})$,
so $\log_{q}\norm{\Ad(t^{2m})}_{\mathcal{L}}\ge2mD\ge2n-2D$. 
\end{remark}

\section{Constructible functions and definable thresholds}

\label{sec:constructible}

In the following let $F$ range over the set $\mathrm{Loc}_{0}$ of
all non-Archimedean local fields of characteristic zero, with no restriction
on the residual characteristic or the ramification.

\subsection{Local fields and the generalized Denef--Pas language }

\label{subsec:language}

We recall the framework of \cite[Sections~2.1--2.5]{MR3732692} and
\cite[Section~1.2]{MR3859937}, to which we refer for more details.
Since we work with families of objects uniformly over $\mathrm{Loc}_{0}$,
we keep track of the subscript $F$ in each notation. In particular,
we fix a uniformizer $\varpi_{F}$ and write $\OO_{F}$ for the ring
of integers of $F$, $\mathfrak{p}_{F}$ for its maximal ideal, $q_{F}$
for the cardinality of the residue field, and $\val_{F}\colon F\to\Z\cup\{\infty\}$
for the valuation, normalized by $\val_{F}(\varpi_{F})=1$. For an
integer $N\ge1$, let $\OO_{F}/N\mathfrak{p}_{F}$ be the quotient
of $\OO_{F}$ by the ideal $N\mathfrak{p}_{F}=\mathfrak{p}^{1+\val_{F}(N)}_{F}$,
so that $\OO_{F}/N\mathfrak{p}_{F}$ is the residue field unless the
residual characteristic divides $N$. Let $\res_{N,F}\colon\OO_{F}\to\OO_{F}/N\mathfrak{p}_{F}$
be the projection and let $\ac_{N,F}\colon F\to\OO_{F}/N\mathfrak{p}_{F}$
be the \emph{angular component map} 
\[
\ac_{N,F}(x)=\res_{N,F}(\varpi^{-\val_{F}(x)}_{F}x)\quad(x\ne0),\,\,\,\,\ac_{N,F}(0)=0.
\]
It depends on the choice of $\varpi_{F}$, and conversely the maps
$\ac_{N,F}$ determine $\varpi_{F}$ \cite[Section~2.4]{MR3732692}.

The \emph{generalized Denef--Pas language} $\gDP$ is a many-sorted
first-order language with the following sorts and symbols \cite[Section~2.2]{MR3732692},
\cite[Definition~1.2.3]{MR3859937}: a sort $\VF$ (the valued field)
and sorts $\RF N$, $N\ge1$ (the residue rings), each carrying a
copy of the language $+,-,\cdot,0,1$ of rings; a sort $\VG\cup\{\infty\}$
(the value group) carrying the language $+,-,0,<$ of ordered abelian
groups, together with a constant $\infty$, and function symbols 
\[
\val\colon\VF\to\VG\cup\{\infty\},\,\,\,\,\ac_{N}\colon\VF\to\RF N\,\,\,(N\ge1).
\]
A local field $F$ with a uniformizer is a structure for $\gDP$:
$\VF$ is interpreted as $F$, $\VG$ as $\Z$, $\RF N$ as $\OO_{F}/N\mathfrak{p}_{F}$,
and the function symbols as above. We do not use further constant
symbols, and in particular no parameters from $F$. Note, however,
that the integers are terms of the sort $\VF$ (sums of $\pm1$),
so that for a fixed non-zero integer $N$ the quantity $\val_{F}(N)$
is definable. 
\begin{definition}[{Definable sets and functions, \cite[Section~2.5]{MR3732692}, \cite[Definition~1.2.4]{MR3859937}}]
 \label{def:definable}~
\begin{enumerate}
\item \label{item:definable-set} A \emph{definable set} is a collection
$X=(X_{F})_{F\in\mathrm{Loc}_{0}}$ (or simply $(X_{F})_{F}$) of
subsets
\[
X_{F}\subseteq F^{m}\times\prod^{\ell}_{t=1}\OO_{F}/N_{t}\mathfrak{p}_{F}\times\Z^{r},\,\,\,F\in\mathrm{Loc}_{0}
\]
with fixed $m,r,\ell\ge0$ and $N_{t}\ge1$, for which there exists an
$\mathcal{L}_{\mathrm{gDP}}$-formula $\varphi$ such that $X_{F}$
is the set defined by $\varphi$ in $F$, for every $F\in\mathrm{Loc}_{0}$. 
\item \label{item:definable-map} If $X$ and $Y$ are definable sets, a
collection $f=(f_{F})_{F}$ of maps $f_{F}\colon X_{F}\to Y_{F}$
is a \emph{definable function} if the collection of the graphs of
the $f_{F}$ is a definable set. 
\end{enumerate}
Since the definable set $\mathrm{VG}$ is simply $\left(\Z\right)_{F}$,
we use the notation $\Z$ instead, so a definable $f:X\rightarrow\VG$
is denoted $f:X\rightarrow\Z$. Moreover, following \cite{MR3859937},
we use the notation $W\subseteq\mathrm{RF}^{*}_{*}$ to mean that
$W\subseteq\prod^{\ell}_{i=1}\mathrm{RF}_{N_{i}}$ for suitable $N_{i}$. 
\end{definition}

Here are some examples. The sets $\Z$, $\N=\{r\in\Z:r\ge0\}$, and
the congruence class $I_{e,M}:=\{r\in\Z:r\equiv e\pmod M\}$ for fixed
integers $M\ge1$ and $e$, are definable subsets, and the function
$f:I_{e,M}\rightarrow\Z$, $r\mapsto(r-e)/M$ is definable. The collection
$\left(\{(x,n)\in F\times\Z:\val_{F}(x)\ge n\}\right)_{F}$ is a definable
subset of $\VF\times\Z$. If $\mathscr{G}$ is a reductive group scheme
over $\bS$, with closed embeddings $\bS\hookrightarrow\A^{D}_{\Z}$
and $\mathscr{G}\hookrightarrow\GL_{d,\bS}\hookrightarrow\A^{d^{2}+1+D}_{\Z}$,
as in \cref{thm:families-intro}, put $S_{F}=\bS(F)$ and $S=(S_{F})_{F\in\mathrm{Loc}_{0}}$.
The collection $\left(\left\{ (g,s):g\in\mathscr{G}_{s}(F),\,s\in S_{F}\right\} \right)_{F}$
is a definable subset of $\VF^{d^{2}+1+D}$. The matrix congruence
subgroups $\matcong n{\rho_{s}}$, the lengths $\ell_{\rho_{s}}$,
and the balls $\matball{\rho_{s}}(r)$ for $r\in\Z_{\ge0}$ are definable
by valuation conditions on the entries of $\rho_{s}(g)$ and its inverse.
The center $Z_{s}$ is definable by commutation with every element
of $G_{s}$. The condition $\ell_{\rho_{s},Z_{s}}(g)=r$ says that
some $z\in Z_{s}$ satisfies $\ell_{\rho_{s}}(zg)=r$ and every $z'\in Z_{s}$
satisfies $\ell_{\rho_{s}}(z'g)\ge r$. Thus $\ell_{\rho_{s},Z_{s}}$
is definable. 

We next introduce constructible functions: 
\begin{definition}[{Constructible functions, \cite[Section~3.1]{MR3732692}, \cite[Definition~1.2.5]{MR3859937}}]
 \label{def:constructible}Let $X$ be a definable set. A collection
$h=(h_{F})_{F}$ of functions $h_{F}\colon X_{F}\to\R$ is \emph{constructible}
if it is a finite sum of functions of the form
\begin{equation}
x\longmapsto \#Z_{F,x}\,q_{F}^{\alpha_{F}(x)}
\prod_{j=1}^{b}\beta_{j,F}(x)
\prod_{i=1}^{d}\frac{1}{1-q_{F}^{c_{i}}},
\qquad F\in\mathrm{Loc}_{0},\quad x\in X_{F}.
\label{eq:constructible}
\end{equation}
Here $\alpha,\beta_{j}\colon X\to\Z$ are definable functions, the
$c_{i}$ are nonzero integers, and $Z_{F,x}$ is the fiber over $x$ of
$Z_{F}$ for a definable set $Z\subseteq X\times\prod_{t=1}^{\ell}\RF{N_{t}}$,
with positive integers $N_{t}$. The nonnegative integers $b,d,\ell$ and the tuple
$(N_{1},\dots,N_{\ell})$ may vary from summand to summand. The constructible functions
on $X$ form a ring $\cO(X)$.
\end{definition}

Definable functions $\beta:X\rightarrow\Z$, the functions $q^{\alpha}_{F}$
with $\alpha:X\rightarrow\Z$ definable, and the characteristic functions
of definable sets are constructible, and the restriction of a constructible
function to a definable subset, extended by zero, is constructible.
The functions $\phi_{a,b,i}$ of \eqref{eq:exponential-polynomials}
below are typical examples. In \cite{MR3859937} most results are
stated for the larger class of functions of $\cO^{\mathrm{exp}}$-class,
which involve an additive character of $F$; no additive character
enters our test functions, and by \cite[Remark~1.2.10]{MR3859937}
the results from \cite{MR3859937} used in this section hold verbatim
for the class $\cO$. 
\begin{definition}[{Constructible conditions, \cite[Definition~1.2.8 and Notation~1.2.9]{MR3859937}}]
 \label{def:condition}Let $X$ be a definable set. A \emph{condition}
on $X$ is a collection $P=(P_{F})_{F}$ of predicates $P_{F}$ on
$X_{F}$. It is \emph{definable} if the collection of the sets $\{x\in X_{F}:P_{F}(x)\text{ holds}\}$
is a definable set, and it is a \emph{constructible condition} (or
\emph{$\mathcal{C}$-condition}), if there is a constructible function
$H$ on $X$ such that 
\[
P_{F}(x)\text{ holds}\iff H_{F}(x)=0,\,\,\text{for every }F\in\mathrm{Loc}_{0},\text{\,\,and\,every\,\,}x\in X_{F}.
\]
\end{definition}

Every definable condition is constructible. However, the converse
fails: the condition $q^{a}_{F}=b$ on $(a,b)\in\Z^{2}$ is the zero
condition of the constructible function $q^{a}_{F}-b$, but it is
not definable \cite[after Definition~1.2.8]{MR3859937}. Constructible
conditions are closed under conjunction, disjunction and universal
quantification over definable sets, but not under negation or existential
quantification \cite[Section~1.3]{MR3859937}. We collect this in
\cref{prop:conditions}, along with a third property which will be
a crucial ingredient in the proof of \cref{thm:main}. 
\begin{proposition}[{Constructible conditions and eventual behavior, \cite[Propositions~1.3.1,1.4.1 and Remark~1.2.10]{MR3859937}}]
 \label{prop:conditions}Let $X,Y$ be definable sets. 
\begin{enumerate}
\item \label{item:condition-positive-boolean} If $P(x)$ and $Q(x)$ are
constructible conditions on $X$, then so are $P(x)\wedge Q(x)$ and
$P(x)\vee Q(x)$. 
\item \label{item:condition-universal} If $P(x,y)$ is a constructible
condition on $X\times Y$, then $\forall y\in Y\colon P(x,y)$ is
a constructible condition on $X$. 
\item \label{item:condition-eventual} (eventual behavior) Let $P(x,r)$
be a constructible condition on $X\times\Z$. Then: 
\begin{enumerate}
\item \label{item:eventual-locus}The condition $Q(x)$ defined by 
\[
\exists r_{0}\in\Z:\forall r\geq r_{0}\,\,\,P(x,r)\text{ holds},
\]
is a constructible condition. Explicitly, for $F\in\mathrm{Loc}_{0}$
and $x\in X_{F}$, the condition $Q_{F}(x)$ means $\exists r_{0}=r_{0}(F,x)$
such that for $r\geq r_{0}$, $P_{F}(x,r)$ holds. 
\item \label{item:eventual-witness} There is a definable function $R:X\rightarrow\Z$
such that for every $F\in\mathrm{Loc}_{0}$, whenever $P_{F}(x,r)$
holds for all sufficiently large $r$, it holds for every $r\ge R_{F}(x)$. 
\end{enumerate}
\end{enumerate}
\end{proposition}

A particularly important property of constructible functions is their
stability under integration. First, we write $\mu_{F}$ for the Haar
measure on $F$, normalized so that $\mu_{F}(\mathcal{O}_{F})=1$.
We use counting measure on $\Z$ and on each finite residue ring $\OO_{F}/N\mathfrak{p}_{F}$.
A definable subset $Y\subseteq\VF^{m}\times\prod^{\ell}_{i=1}\RF{N_{i}}\times\Z^{r}$
inherits a family of measures $\mu_{Y,F}$ on $Y_{F}$ by restriction
from the corresponding product measure on the ambient space $F^{m}\times\prod^{\ell}_{i=1}\OO_{F}/N_{i}\mathfrak{p}_{F}\times\Z^{r}$.
In particular, $\mu_{Y,F}$ vanishes if $Y_{F}$ has empty interior in the
ambient product.
For $f\in\mathcal{C}(Y)$, integration means ordinary Lebesgue integration
of each specialization $f_{F}$. For integrable $f_{F}$, writing
$\widetilde{f}_{F}$ for its extension by zero to the ambient product,
we have 
\[
\int_{Y_{F}}f_{F}\,d\mu_{Y,F}:=\sum_{z\in\Z^{r}}\sum_{\xi\in\prod^{\ell}_{i=1}\OO_{F}/N_{i}\mathfrak{p}_{F}}\int_{F^{m}}\widetilde{f}_{F}(u,\xi,z)\,du.
\]

\begin{theorem}[{Integration preserves constructibility, \cite[Theorem~4.1.1]{MR3732692}}]
 \label{thm:integration}For definable $X,Y$ and $W\subseteq X\times Y$,
set $W_{F,x}=\{y\in Y_{F}:(x,y)\in W_{F}\}$. Then for every $f\in\mathcal{C}(W)$,
there exists $g\in\mathcal{C}(X)$ such that, for every $F\in\mathrm{Loc}_{0}$
and $x\in X_{F}$, 
\[
\int_{W_{F,x}}\left|f_{F}(x,y)\right|\,d\mu_{Y,F}(y)<\infty\,\,\,\Longrightarrow\,\,g_{F}(x)=\int_{W_{F,x}}f_{F}(x,y)\,d\mu_{Y,F}(y).
\]
In particular, if every fiber is integrable, then the family of fiber
integrals itself belongs to $\mathcal{C}(X)$. 
\end{theorem}

\subsection{The model-theoretic statement}

The definable structure of constructible functions makes it possible
to control certain analytic properties uniformly as the local field
varies. Here we are concerned with \emph{eventual vanishing.} Suppose
$h\in\mathcal{C}(X\times\Z_{\geq0})$ is a constructible function,
and suppose that for every $F\in\mathrm{Loc}_{0}$ and $x\in X_{F}$,
$h_{F}(x,\cdot)$ vanishes when $r\gg_{F,x}1$. It follows from \cref{prop:conditions}
that there is a \emph{definable threshold} for the vanishing of $h$,
namely, there is a definable function $R:X\rightarrow\Z$ such that
$h_{F}(x,r)=0$ for $r>R_{F}(x)$, uniformly in $F\in\mathrm{Loc}_{0}$
and $x\in X_{F}$. The basic mechanism behind this phenomenon is a
preparation result explained below: after a definable reparameterization
by residue-ring variables, the domain of any constructible $\widetilde{h}\in\mathcal{C}(X\times\Z_{\geq0})$
can be partitioned into cells whose fibers in the $r$-variable are
either finite, with definable endpoints and hence definably bounded
diameter, or unbounded. On each fiber the reparameterized function
is a linear combination of exponential-polynomials on a congruence
class. Such a combination vanishes eventually on an unbounded fiber
if and only if it vanishes identically there. Taking $\widetilde{h}=h$
therefore gives a definable upper bound for the support of $h_{F}(x,\cdot)$.

For the moment suppress the parameter $s$, writing $\rho$ for $\rho_{s}$.
As explained in the introduction, bounding the support of $\matcong n{\rho}$-invariant
matrix coefficients of supercuspidal representations can be reduced
to analyzing the support of complex linear combinations of orbital
integrals $O_{y}(\shelltest_{n,r})$, where $y$ is in $Y_{F}$, the
set of regular semisimple compact elements of $\bG(F)$. In \cref{sec:orbital}
we show that $u(y,n,r)=O_{y}(\shelltest_{n,r})$ is a constructible
function on $Y\times\Z_{\geq1}\times\Z_{\geq0}$. This gives a direct
link to the question of eventual vanishing of constructible functions.
However, there are two obstacles: 
\begin{enumerate}
\item \label{item:obstacle-complex-span} The orbital integrals themselves
$O_{y}(\shelltest_{n,r})$ \emph{need not eventually vanish} in $r$,
only a suitable complex linear combination $\sum^{m}_{i=1}a_{F,n,i}O_{y_{i}}(\shelltest_{n,r})$
of orbital integrals, with unspecified coefficients. While the resulting
family need not be constructible, it is still amenable to the techniques
of this section, as such linear combinations are still piecewise exponential
polynomials on congruence classes. 
\item \label{item:obstacle-orbital-parameter} For $h\in\mathcal{C}(Y\times\Z_{\geq1}\times\Z_{\geq0})$,
eventual vanishing of $h_{F}(y,n,\cdot)$ for each $y$ does not guarantee
a definable threshold independent of $y$. Indeed, for $y\ne0$, $h(y,n,r):=1_{\left\{ r\leq\val_{F}(y)+n\right\} }$
eventually vanishes as $r\gg_{F,y,n}1$, but its support endpoints
tend to infinity as $\val_{F}(y)$ tends to infinity. 
\end{enumerate}
By Clozel's theorem \cite{MR986793}, for each fixed $F$ and $n$, the sequences
$r\mapsto u_{F}(y,n,r)$, with $y\in Y_{F}$, span a finite-dimensional
space. This allows us, for each $F$ and $n$, to
reduce to finitely many $y$'s, and produce a threshold independent
of $y$, but a priori not necessarily definable. The final ingredient
is to reduce to a constructible condition, and produce a definable
threshold from \cref{prop:conditions}.

With the above discussion in mind, here is the main proposition of
this section: 
\begin{proposition}[Definable threshold]
\label{prop:threshold} Let $Y,Z$ and $X\subseteq Y\times Z$ be
definable sets, and let $u\in\mathcal{C}\bigl(X\times\Z_{\ge0}\bigr)$.
For $z\in Z_{F}$, put $Y_{F,z}:=\{y\in Y_{F}:(y,z)\in X_{F}\}$.
Suppose that for every $F\in\mathrm{Loc}_{0}$ and $z\in Z_{F}$,
the space 
\[
W_{F,z}:=\operatorname{span}_{\C}\left\{ (u_{F}(y,z,r))_{r\geq0}:y\in Y_{F,z}\right\} 
\]
is finite-dimensional. Then there is a definable function $R:Z\to\Z_{\geq0}$
such that, for every $F\in\mathrm{Loc}_{0}$, $z\in Z_{F}$, every
$a\in W_{F,z}$, 
\[
a(r)=0\text{ for }r\gg1\Longrightarrow\,\,a(r)=0\text{ for every }r\geq R_{F}(z).
\]
Moreover, if $Z=S\times\Z_{\ge1}$, for a definable set $S$, there
is an integer $A\geq0$, and a definable function $\gamma:S\rightarrow\Z$,
such that for every $F\in\mathrm{Loc}_{0}$, $s\in S_{F}$ and $n\in\Z_{\geq1}$:
\[
R_{F}(s,n)\leq An+\gamma_{F}(s).
\]
\end{proposition}

\begin{remark}
\label{rem:dimension-free}No bound on $\dim W_{F,z}$ uniform in
$F,z$ is assumed, and no definable choice of a basis of $W_{F,z}$
is made. A basis is used for a fixed $F,z$ only to establish eventual
membership in the prepared tail space. 
\end{remark}

\begin{remark}[Families of groups]
For the family in \cref{thm:families-intro}, we apply \cref{prop:threshold}
with $z=(s,n)$, where $s\in\bS(F)$ specifies the group $G_{s}$ and $n$
the congruence level. For fixed $F,s,n$, the element $y$ ranges over the
compact regular semisimple elements of $G_{s}$. The threshold $R_{F}(s,n)$
is common to all eventually vanishing linear combinations of the sequences
$r\mapsto O_{y}(\shelltest_{s,n,r})$. Once the fiber is fixed, this
threshold depends only on $n$.
\end{remark}

Moreover, with the above discussion in mind, the following lemma,
due to Cluckers--Gordon--Halupczok, explains where the linearity
in \cref{thm:main} comes from. We include a direct proof for the
convenience of the reader. 
\begin{lemma}[Definable functions are affinely bounded]
 \label{lem:vg-bound}Let $Z=S\times\Z_{\ge1}$, for a definable
set $S$, and let $R:Z\rightarrow\Z$ be a definable function. Then
there is $A\in\Z_{\geq0}$, and a definable function $\gamma:S\rightarrow\Z$,
such that for every $F\in\mathrm{Loc}_{0}$, $s\in S_{F}$ and $n\in\Z_{\geq1}$:
\[
\left|R_{F}(s,n)\right|\leq An+\gamma_{F}(s).
\]
Moreover, if $Z=\Z_{\ge1}$ ($S$ is a singleton), there are integers
$A,B\geq0$, and $N\geq1$, such that for every $F\in\mathrm{Loc}_{0}$
and $n\in\Z_{\geq1}$: 
\[
\left|R_{F}(n)\right|\leq An+B+\val_{F}(N).
\]
\end{lemma}

\begin{proof}
By \cite[Definition~5.2.2 and Corollary~5.2.5, pp.~72--75]{MR3732692}
(see also \cite[Proof of Theorem~2.1.1, pp.~139--140]{MR3859937}),
there is a finite definable partition $Z=\bigsqcup^{t}_{i=1}P_{i}$,
such that $R|_{P_{i}}$ is ``approximately linear'' in $n$ in the
following sense. There are $M\in\Z_{\geq1}$, $a_{i},b_{i}\in\Q$,
and definable functions $\gamma_{i}\colon S\to\Z$, such that for
every $F\in\mathrm{Loc}_{0}$ and every $(s,n)\in P_{i,F}$: 
\[
0\le R_{F}(s,n)-a_{i}\bigl(\gamma_{i,F}(s)+b_{i}n\bigr)\le\val_{F}(M).
\]
Take integers $A\ge\max_{1\le i\le t}\left|a_{i}b_{i}\right|$ and
$D\ge\max_{1\le i\le t}\left|a_{i}\right|$ and consider: 
\[
\gamma_{F}(s):=D\max_{1\le i\le t}\left|\gamma_{i,F}(s)\right|+\val_{F}(M),\,\,\,\text{ for }F\in\mathrm{Loc}_{0}.
\]
Since $\gamma\colon S\to\Z_{\ge0}$ is definable, this gives the first
assertion. If $S$ is a singleton, then $d_{F}:=\gamma_{F}(s)$ is
a definable integer, depending only on $F$. By the paragraph following
\cite[Theorem~2.1.1, p.~139]{MR3859937}, there are integers $B\ge0$
and $N\ge1$, such that for every $F\in\mathrm{Loc}_{0}$: 
\[
d_{F}\le B+\val_{F}(N).
\]
This proves the second assertion. 
\end{proof}

\subsubsection{Proof of \cref{prop:threshold}}

Let $u(x,r)\in\mathcal{C}\bigl(X\times\Z_{\ge0}\bigr)$. We first
apply a preparation result \cite[Proposition~1.4.2]{MR3859937}, which
allows us to write $u$, after a residue-ring reparameterization and
a finite definable partition of the resulting domain, as an exponential
polynomial in $r$ on each cell (\eqref{eq:cell-preparation} below).
Concretely, by \cite[Proposition~1.4.2]{MR3859937}, there exists
a definable bijection: 
\[
\sigma:X\times\Z_{\geq0}\longrightarrow U^{\mathrm{par}}\subseteq\mathrm{RF}^{*}_{*}\times X\times\Z_{\ge0},\,\,\,\,\,\sigma(x,r)=(\xi(x,r),x,r),
\]
and a finite definable partition into cells: 
\[
U^{\mathrm{par}}=\bigsqcup^{J}_{j=1}A_{j}\,\,\,\,,\,\,\,\,A_{j}=\{(\xi,x,r):(\xi,x)\in B_{j},\,r\in A_{j,\xi,x}\},
\]
where 
\begin{equation}
A_{j,\xi,x}:=\{r\in\Z_{\geq0}:\alpha_{j}(\xi,x)\leq r\leq\beta_{j}(\xi,x),\ r\equiv e_{j}\pmod{m_{j}}\}.\label{eq:presburger-cell}
\end{equation}
Here, $m_{j}\ge1,e_{j}$ are fixed integers. $B_{j}$ is the projection
of $A_{j}$ to $\mathrm{RF}^{*}_{*}\times X$, $\alpha_{j}:B_{j}\rightarrow\mathrm{VG}$
is definable, and $\beta_{j}$ is either definable, or $\beta_{j}=\infty$.
On each cell $A_{j}$, the function $u\circ\sigma^{-1}$ has the following
form:\label{On-each-cells} 
\begin{equation}
(u\circ\sigma^{-1})_{F}|_{A_{j,F}}(\xi,x,r)=\sum^{N_{j}}_{i=1}c_{ij,F}(\xi,x)r^{a_{i,j}}q^{b_{i,j}r}_{F},\label{eq:cell-preparation}
\end{equation}
where $c_{ij,F}(\xi,x)$ are coefficients independent of $r$, $b_{i,j}\in\Q$
and integers $0\leq a_{i,j}\leq d_{0}$, for a fixed integer $d_{0}\ge0$.

Let $M\geq1$ be divisible by all $m_{j}$ in \eqref{eq:presburger-cell},
as well as the denominators of all $b_{ij}$'s in \eqref{eq:cell-preparation}.
Let $\mathcal{B}\subseteq\Q$ be the collection of all $b_{i,j}$.
For $0\le a<M$, $b\in\mathcal{B}$ and $0\leq i\leq d_{0}$, set:
\begin{equation}
\phi_{a,b,i,F}(r):=\begin{cases}
k^{i}q^{Mbk}_{F} & r=Mk+a\text{ for }k\in\Z_{\geq0}\\
0 & \text{otherwise.}
\end{cases}.\label{eq:exponential-polynomials}
\end{equation}
Note that $Mb$ is an integer, and by \cref{def:constructible}, each
$\phi_{a,b,i}$ belongs to $\mathcal{C}(\Z_{\ge0})$. Denote 
\[
E_{F}:=\operatorname{span}_{\C}\{\phi_{a,b,i,F}\}_{a,b,i},
\]
namely the space of all finite linear combinations $\sum_{a,b,i}\alpha_{a,b,i}\phi_{a,b,i,F}$,
for $\alpha_{a,b,i}\in\C$. 
\begin{lemma}
\label{lem:tail-space}For every $F\in\mathrm{Loc}_{0}$, and every
$x\in X_{F}$, there exists $r_{0}=r_{0}(F,x)\in\N$ such that if
$r>r_{0}$: 
\begin{enumerate}
\item \label{item:tail-infinite-cells} If $\sigma(x,r)=(\xi,x,r)\in A_{j}$,
then $A_{j,\xi,x}$ is infinite. 
\item \label{item:tail-membership} $\widetilde{u}:=u_{F}(x,\cdot)$ agrees
with an element of $E_{F}$ on the set $\left\{ r>r_{0}\right\} $. 
\end{enumerate}
\end{lemma}

\begin{proof}
Fix $F$ and $x\in X_{F}$. The residue rings are finite, so only
finitely many tuples $(\xi,j)$ occur over $x$. We take $r_{0}=\underset{\xi,j}{\max}\left\{ \left\{ \alpha_{j}(\xi,x)\right\} \cup\left\{ \beta_{j}(\xi,x):\beta_{j}<\infty\right\} \right\} +1$,
where $(\xi,j)$ ranges over these tuples. $r_{0}$ is larger than
the endpoint of all finite cell fibers $A_{j,\xi,x}$, so $A_{j,\xi,x}$
is infinite if it contains some $r>r_{0}$. In the latter case, we
can ignore the start $\alpha_{j}(\xi,x)$ of the cell, so that $A_{j,\xi,x}\cap\{r>r_{0}\}$
is a union of congruence classes modulo $M$ restricted to $r>r_{0}$.
By \eqref{eq:cell-preparation}, $\widetilde{u}$ agrees with a linear
combination of the functions in \eqref{eq:exponential-polynomials},
when $r>r_{0}$.
\end{proof}

The next lemma gives a weakened version of \cref{prop:threshold},
with a threshold $r_{1}(F,z)$ independent of $y$ in place of the
definable threshold $R_{F}(z)$. 
\begin{lemma}
\label{lem:uniform-tail}For every $F\in\mathrm{Loc}_{0}$, and every
$z\in Z_{F}$ there is $r_{1}(F,z)\in\N$ such that for every $a\in W_{F,z}$,
one has: 
\begin{enumerate}
\item \label{item:uniform-tail-membership} $a|_{r>r_{1}(F,z)}\in E_{F}|_{r>r_{1}(F,z)}$. 
\item \label{item:uniform-tail-vanishing} $a(r)=0\,\forall r\gg1\Longrightarrow a(r)=0\,\forall r>r_{1}(F,z)$. 
\end{enumerate}
\end{lemma}

\begin{proof}
Let $F\in\mathrm{Loc}_{0}$, let $z\in Z_{F}$, and let $y_{1},...,y_{t}\in Y_{F,z}$
be such that $\left\{ u_{F}(y_{j},z,\cdot)\right\} ^{t}_{j=1}$ is
a basis for $W_{F,z}$. By \cref{lem:tail-space}, choose $e_{j}\in E_{F}$,
such that $e_{j}(r)=u_{F}(y_{j},z,r)$ for $r>r_{0}(F,y_{j},z)$.
Set $r_{1}(F,z):=\max\bigl(\{0\}\cup\{r_{0}(F,y_{j},z):1\leq j\leq t\}\bigr)$.
Hence, for every $\boldsymbol{\alpha}:=(\alpha_{1},\dots,\alpha_{t})\in\C^{t}$,
we have 
\[
w_{\boldsymbol{\alpha}}(r):=\sum^{t}_{j=1}\alpha_{j}\cdot u_{F}(y_{j},z,r)=\sum^{t}_{j=1}\alpha_{j}e_{j}(r),\,\text{ for every }r>r_{1}(F,z).
\]
This proves \cref{item:uniform-tail-membership}. Now suppose that
$w_{\boldsymbol{\alpha}}(r)=0$ for $r\gg1$. In particular,
$\sum^{t}_{j=1}\alpha_{j}e_{j}(r)=0$ for $r\gg1$. Since $E_{F}:=\operatorname{span}_{\C}\{\phi_{a,b,i,F}\}_{a,b,i}$,
on each congruence class $\left\{ kM+a\right\} $, there are suitable
constants $c_{\boldsymbol{\alpha},a,b,i}\in\C$ such that: 
\begin{equation}
\sum^{t}_{j=1}\alpha_{j}e_{j}(kM+a)=\sum_{b,i}c_{\boldsymbol{\alpha},a,b,i}\cdot\phi_{a,b,i,F}(kM+a)=\sum_{b,i}c_{\boldsymbol{\alpha},a,b,i}\cdot k^{i}q^{Mbk}_{F}.\label{eq:tail-polynomial}
\end{equation}
By, for example, \cite[Lemma 2.1.7]{MR3210968} applied separately to
the real and imaginary parts, \eqref{eq:tail-polynomial} vanishes
for $k\gg1$ if and only if $c_{\boldsymbol{\alpha},a,b,i}=0$
for all $b,i$. We deduce that $\sum^{t}_{j=1}\alpha_{j}e_{j}(r)=0$
for every $r\geq0$, and hence $w_{\boldsymbol{\alpha}}(r)=0$
for $r>r_{1}(F,z)$. 
\end{proof}

For each $T\in\Z_{\geq0}$ consider the following condition: 
\[
P_{F}(z,T)\,\,:\,\,\forall y\in Y_{F,z}:u_{F}(y,z,\cdot)|_{r\geq T}\in E_{F}|_{r\geq T}.
\]

\begin{lemma}
\label{lem:determinant-condition}The collection $P(z,T):=\left\{ P_{F}(z,T)\right\} _{F\in\mathrm{Loc}_{0}}$
is a $\mathcal{C}$-condition on $Z\times\Z_{\geq0}$. 
\end{lemma}

\begin{proof}
Enumerate $\left\{ \phi_{a,b,i}\right\} _{a,b,i}$ as $\left\{ \phi_{j}\right\} ^{D}_{j=1}$,
where $D=M\cdot\left|\mathcal{B}\right|\cdot(d_{0}+1)=\dim E_{F}$
is independent of $F$. Fix $F\in\mathrm{Loc}_{0}$, $x\in X_{F}$.
Note that $u_{F}(x,\cdot)|_{r\geq T}\in E_{F}|_{r\geq T}$ $\Longleftrightarrow$
there exists $\alpha_{1},\dots,\alpha_{D}\in\C$ such that $u_{F}(x,r)=\sum^{D}_{j=1}\alpha_{j}\phi_{j,F}(r)$
for $r\geq T$. In particular, denoting 
\[
H_{F}(x,t_{0},\dots,t_{D}):=\det\begin{pmatrix}\phi_{1,F}(t_{0}) & \cdots & \phi_{D,F}(t_{0}) & u_{F}(x,t_{0})\\
\vdots & \, & \vdots & \vdots\\
\phi_{1,F}(t_{D}) & \cdots & \phi_{D,F}(t_{D}) & u_{F}(x,t_{D})
\end{pmatrix},
\]
we get that $u_{F}(x,\cdot)|_{r\geq T}\in E_{F}|_{r\geq T}$ implies
that $H_{F}(x,t_{0},\dots,t_{D})=0$ for every $(t_{0},\dots,t_{D})\in\Z^{D+1}_{\geq T}$.
The converse holds as well. Indeed, by the proof of the previous lemma,
the functions $\phi_{1,F}|_{r\geq T},\dots,\phi_{D,F}|_{r\geq T}$
are linearly independent. If $u_{F}(x,\cdot)|_{r\geq T}\notin E_{F}|_{r\geq T}$
then the functions $\left\{ \phi_{j,F}|_{r\geq T}\right\} ^{D}_{j=1}\cup\left\{ u_{F}(x,\cdot)|_{r\geq T}\right\} $
are linearly independent. In particular, the image of the map $t\mapsto\overrightarrow{v}(t):=\left(\phi_{1,F}(t),\dots,\phi_{D,F}(t),u_{F}(x,t)\right)$
for $t\geq T$, spans $\C^{D+1}$. Therefore, we can find $t_{0},\dots,t_{D}\in\Z_{\geq T}$
such that the vectors $\left\{ \overrightarrow{v}(t_{j})\right\} ^{D}_{j=0}$
are linearly independent. Hence, $H_{F}(x,t_{0},\dots,t_{D})\neq0$.
Extend $H$ by zero outside $X\times\Z^{D+1}_{\geq0}$. We deduce
that: 
\begin{align*}
P_{F}(z,T)\text{ holds } & \Longleftrightarrow\left\{ \forall y\in Y_{F,z},\,\forall t_{0},\dots,t_{D}\in\Z_{\geq T}\,\,H_{F}(y,z,t_{0},\dots,t_{D})=0\right\} \\
 & \Longleftrightarrow\left\{ \begin{aligned} & \forall y\in Y_{F},\,\forall t_{0},\dots,t_{D}\in\Z_{\geq0}\\
 & \quad1_{\left\{ t_{i}\geq T\,\,0\leq i\leq D\right\} }1_{\left\{ y\in Y_{F,z}\right\} }\cdot H_{F}(y,z,t_{0},\dots,t_{D})=0
\end{aligned}
\right\} .
\end{align*}
Since the function $1_{\left\{ t_{i}\geq T\,\,0\leq i\leq D\right\} }1_{\left\{ y\in Y_{F,z}\right\} }\cdot H_{F}(y,z,t_{0},\dots,t_{D})$
is constructible, by \cref{prop:conditions}\cref{item:condition-universal},
$P(z,T)$ is a $\mathcal{C}$-condition, as required. 
\end{proof}

We can now finish the proof of the proposition: 
\begin{proof}[Proof of \cref{prop:threshold}]
Extend $P(z,T)$ to $T\in\Z$ by $P(z,\max\{0,T\})$. Consider the
condition 
\[
Q_{F}(z):\,\exists T_{0}(F,z)\in\Z_{\geq0}\text{ such that }\forall T\geq T_{0},\,P_{F}(z,T)\text{ holds.}
\]
By \cref{lem:determinant-condition}, $P(z,T)$ is a $\mathcal{C}$-condition.
By \cref{prop:conditions}\cref{item:condition-eventual}, it follows
that $Q(z)$ is a $\mathcal{C}$-condition, and there exists a definable
$R:Z\rightarrow\Z_{\geq1}$ such that for every $F\in\mathrm{Loc}_{0}$
and every $z\in Z_{F}$, if $Q_{F}(z)$ holds, then $P_{F}(z,T)$
holds for every $T\geq R_{F}(z)$.

Fix $F$ and $z$. By \cref{lem:uniform-tail}\cref{item:uniform-tail-membership},
taking $T_{0}:=r_{1}(F,z)+1$ we get that $P_{F}(z,T)$ holds for
every $T\geq T_{0}$, so $Q_{F}(z)$ holds. Hence, we deduce by linearity
that for every $a\in W_{F,z}$, $a|_{r\geq R_{F}(z)}\in E_{F}|_{r\geq R_{F}(z)}$.
Suppose now that $a(r)=0$ for $r\gg1$. By the proof of \cref{lem:uniform-tail}\cref{item:uniform-tail-vanishing},
any element of $E_{F}$ vanishes eventually if and only if it vanishes
identically. Thus, $a|_{r\geq R_{F}(z)}=0$, as required.

The second part of the proposition follows by \cref{lem:vg-bound}. 
\end{proof}

\section{Orbital integrals and positive shells}

\label{sec:orbital}

Our goal in this section is to show that the functions $\shelltest_{n,r}$
and the orbital integrals $O_{y}(\shelltest_{n,r})$ introduced in
\cref{sec:intro} are constructible. We work on a finite partition
of the parameter scheme into locally closed pieces on which the fiber
dimension and absolute rank are constant. The constructible functions
on the pieces then give a constructible function on the whole family.
We first choose invariant measures on the orbits so that $u(y,n,r):=O_{y}(\shelltest_{n,r})$
is constructible.

Constructibility of orbital integrals has been established in several
(group-theoretic, and model-theoretic) settings. Cluckers--Gordon--Halupczok
treat Lie-algebra orbital integrals in \cite[Proposition 4.3 and Corollary 4.4]{MR3297157},
and semisimple group orbital integrals in \cite[Appendix B, Lemma 14.15]{MR3437869}.
These results are formulated in sufficiently large residue characteristic,
due to the model-theoretic framework. In \cite{MR3859937}, the model-theoretic
framework works in $\mathrm{Loc}_{0}$, but the orbital-integral statement
\cite[Lemma 4.3.2 and Remark 4.3.3]{MR3859937} still imposes a lower
bound on the residue characteristic and concerns Lie-algebra orbital
integrals.

We are not aware of a result in the literature covering precisely
our setting. We therefore prove in \cref{lem:orbital-density} below
that for a reductive group scheme over an affine finite-type $\Z$-scheme,
regular semisimple group orbits admit invariant measures whose integrals
form constructible families over all fields in $\Loc$.

\subsection{Canonical measures}

\label{subsec:Canonical-measures-on}

We start with an example which illustrates the potential normalization
issue. 
\begin{example}
\label{exa:Weils canonical measure}Consider the scheme $X=\mathrm{Spec}\left(\Z[x,y,(x-1)^{-1}]/(y^{2}-x^{3}+x)\right)$.
This is a smooth affine $\Z$-scheme of finite type. Its generic fiber is obtained
by removing $(1,0)$ and $\infty$ from the elliptic curve $y^2=x^3-x$.
The family $(X(\OO_F))_F$ is $\gDP$-definable. On each $X(\OO_F)$ we obtain
a probability measure $d\mu^{\mathrm{naive}}_{X(\OO_F)}$ as the inverse limit
of the uniform probability measures on $X(\OO_F/\mathfrak p_F^k)$
(e.g.~for $X=\mathrm{SL}_{n}$, this would simply mean taking the
normalized Haar probability measure on $\mathrm{SL}_{n}(\mathcal{O}_{F})$).
Now, take $h_{F}:=1_{\left\{ (x,y)\in X(\mathcal{O}_{F}):\val_{F}(x),\val_{F}(y)\geq1\right\} }$.
Then $h$ is a constructible function on $\left(X(\mathcal{O}_{F})\right)_{F}$,
but on the other hand, Hensel's lemma implies that: 
\begin{equation}
d\mu^{\mathrm{naive}}_{X(\mathcal{O}_{F})}(h_{F})=\frac{1}{\#X(k_{F})},\text{ for every }F\in\mathrm{Loc}_{0}.\label{eq:non motivic measure}
\end{equation}
Here $k_{F}=\OO_{F}/\mathfrak{p}_{F}$ is the residue field.
This normalization introduces the \textbf{reciprocal} of a point count,
and \cref{def:constructible} does not guarantee that such reciprocals
are constructible, so stability under integration with these probability
measures is not guaranteed.
\end{example}

To obtain a normalization compatible with constructible integration,
we use definable volume forms. We recall the Serre--Oesterlé measure
associated with an embedded $F$-analytic submanifold.
\begin{definition}[{Serre--Oesterlé canonical measure, \cite{MR644559,MR656627} and \cite[Section 8.2]{MR2630060}}]
 \label{def:Serre--Oesterle}Let $F$ be a non-archimedean local
field, and let $X\subseteq F^{N}$ be an $m$-dimensional $F$-analytic
submanifold. For each $I=\{i_{1}<\cdots<i_{m}\}\subseteq\{1,\ldots,N\}$,
write $\omega_{I}=(dx_{i_{1}}\wedge\cdots\wedge dx_{i_{m}})|_{X}$.
Fix an ordering of these subsets. For each $I$, set $X_{I}:=\{x\in X:\omega_{I}(x)\ne0\}$,
so that $\omega_{I}|_{X_{I}}$ is invertible, and define analytic
functions $a_{J,I}:X_{I}\to F$ by $\omega_{J}|_{X_{I}}=a_{J,I}\omega_{I}|_{X_{I}}$.
We set: 
\[
U_{I}=\left\{ x\in X_{I}:\begin{array}{ll}
\left|a_{J,I}(x)\right|_{F}<1 & \text{for every }J<I,\\
\left|a_{J,I}(x)\right|_{F}\leq1 & \text{for every }J>I
\end{array}\right\} .
\]
Since at every point $x\in X$ at least one of the forms $\omega_{I}(x)$
is nonzero, these sets form a finite partition of $X$. Each~$U_{I}$
is open, since the ratios are analytic, and hence it is also closed.
We obtain an invertible $F$-analytic top form $\omega_{X}$ on $X$
by setting $\omega_{X}|_{U_{I}}:=\omega_{I}|_{U_{I}}$. The associated
\emph{Serre--Oesterlé canonical measure} satisfies 
\[
d\mu^{\mathrm{can}}_{X}=\left|\omega_{X}\right|_{F}=\max_{\left|I\right|=m}\left|\omega_{I}\right|_{F}.
\]
The measure $d\mu^{\mathrm{can}}_{X}$ depends on the embedding $X\subseteq F^{N}$,
but is independent of the chosen ordering.
\end{definition}

The Serre--Oesterlé measures are compatible with the theory of motivic
integration. For a definable family $\left\{ X_{z}\right\} _{z}$
of such embedded manifolds, this rule gives a definable family $\left\{ \omega_{X_{z}}\right\} _{z}$
of top forms, and measures $\left\{ d\mu^{\mathrm{can}}_{X_{z}}\right\} _{z}$
and \cref{thm:integration} extends to the following (see e.g.~\cite[discussion after Proposition 4.3.2]{MR3732692}
and \cite[Proof of Proposition-Definition 3.3.1]{MR3997895}). 
\begin{proposition}[Integration against canonical measures]
 \label{prop:canonical-integration}Let $Z$ be an $\gDP$-definable
set, let $W\subseteq\mathrm{VF}^{n}\times Z$ be an $\gDP$-definable
subset, and let $m\in\N$. For $F\in\operatorname{Loc}_{0}$ and $z\in Z_{F}$,
assume that $X_{F,z}:=\{x\in F^{n}:(x,z)\in W_{F}\}$ is an $m$-dimensional
submanifold of $F^{n}$, and let $d\mu^{\mathrm{can}}_{X_{F,z}}$
be the corresponding canonical measure. Then for every $f\in\mathcal{C}(W)$,
there exists $g\in\mathcal{C}(Z)$ such that 
\[
g_{F}(z)=\int_{X_{F,z}}f_{F}(x,z)\,d\mu^{\mathrm{can}}_{X_{F,z}}(x).
\]
for every $F\in\operatorname{Loc}_{0}$ and $z\in Z_{F}$ for which
$f_{F}(\cdot,z)$ is absolutely integrable. In particular, if every
fiber restriction is integrable, this identity holds for all $F$
and $z$. 
\end{proposition}

\subsection{Orbital integrals}

\label{subsec:Orbital-integrals}

Let $\bS$ be an affine scheme of finite type over $\Z$, let $\psi:\bH\to\bS$
be a reductive group scheme. Write $S=(\bS(F))_{F\in\mathrm{Loc}_{0}}$,
$H=(\bH(F))_{F\in\mathrm{Loc}_{0}}$ and $H^{\mathrm{rss}}=(H^{\mathrm{rss}}_{F})_{F\in\mathrm{Loc}_{0}}$,
where 
\[
H^{\mathrm{rss}}_{F}:=\{y\in H_{F}:y\in\bH_{\psi_{F}(y)}(F)^{\mathrm{rss}}\}.
\]
Note that $S,H$ and $H^{\mathrm{rss}}$ are $\gDP$-definable sets.
For $y\in H^{\mathrm{rss}}_{F}$, put $s=\psi_{F}(y)$ and set 
\[
\Orbit_{y}:=\left\{ gyg^{-1}:g\in\bH_{s}(F)\right\} ,
\]
the corresponding $\bH_{s}(F)$ orbit. For every $F\in\mathrm{Loc}_{0}$,
each $\Orbit_{y}$ is an $F$-analytic submanifold. 
\begin{lemma}[Constructible orbital integrals]
 \label{lem:orbital-density}Let $\psi:\bH\to\bS$ be a morphism
as above, and write $s=\psi_{F}(y)$ for $y\in H^{\mathrm{rss}}_{F}$.
There exists a collection $\{d\mu_{y}\}_{F,y\in H^{\mathrm{rss}}_{F}}$
of nonzero positive, $\bH_{s}(F)$-invariant measures on the conjugacy
orbits $\Orbit_{y}$ with the following property. Let $A$ be a definable
set and $f\in\mathcal{C}(H\times A)$, and assume that, for every
$F\in\mathrm{Loc}_{0}$, $s\in\bS(F)$ and $a\in A_{F}$, the restriction
of $f_{F}(\cdot,a)$ to $\bH_{s}(F)$ is compactly supported and locally
constant. Then there exists $h\in\mathcal{C}(H^{\mathrm{rss}}\times A)$
such that, for every $F\in\mathrm{Loc}_{0}$ and $(y,a)\in H^{\mathrm{rss}}_{F}\times A_{F}$,
\[
h_{F}(y,a)=O_{y}(f_{F}(\cdot,a))=\int_{\Orbit_{y}}f_{F}(x,a)\,d\mu_{y}(x).
\]
\end{lemma}

The normalization of $d\mu_{y}$ depends on $y$, but this will play
no role. 
\begin{proof}
Fix an affine embedding $\bH\hookrightarrow\mathbb{A}^{N}_{\bS}$
over $\bS$, and consider 
\[
W_{F}=\{(x,y)\in H_{F}\times_{S_{F}}H^{\mathrm{rss}}_{F}:\exists g\in\bH_{\psi_{F}(y)}(F),\ x=gyg^{-1}\}.
\]
The set $W$ is definable. For $y\in H^{\mathrm{rss}}_{F}$, put $s=\psi_{F}(y)$.
Note that $\bH_{s}$ is a reductive group over a field $F$ of characteristic
zero, so it is smooth. The fibers of $W\to H^{\mathrm{rss}}$ are
the orbits $\Orbit_{y}$, which are $F$-analytic submanifolds of
$F^{N}$. After a finite definable partition of $H^{\mathrm{rss}}$
according to $\dim\Orbit_{y}=\operatorname{rank}_{F}\bigl(\Ad(y)-\mathrm{id}_{\Lie(\bH_{s})(F)}\bigr)$,
we may assume that their dimension $m$ is constant.

\textbf{Step 1}: We choose uniformly definable analytic invertible
top forms $\omega_{\Orbit_{y}}$ on $\Orbit_{y}$ such that $\left\{ d\mu^{\mathrm{can}}_{\Orbit_{y}}\right\} _{y,F}=\left\{ \left|\omega_{\Orbit_{y}}\right|_{F}\right\} _{y,F}$
is the family of canonical measures on $\Orbit_{y}$ defined in \cref{def:Serre--Oesterle}.
These measures satisfy that $\int_{\Orbit_{y}}f_{F}(x,a)\cdot d\mu^{\mathrm{can}}_{\Orbit_{y}}(x)$
is constructible, but they need not be $\bH_{s}(F)$-invariant.

\textbf{Step 2}: We produce $\bH_{s}(F)$-invariant measures on these
orbits. Let $\mathbf{C}_{y}:=C_{\bH_{s}}(y)$ be the algebraic centralizer,
and set $\alpha_{y}:=\omega_{\Orbit_{y}}|_{y}\in\bigwedge^{m}T^{*}_{y}\Orbit_{y}$.
The torus $\mathbf{C}^{\circ}_{y}$ acts trivially on this one-dimensional
space, since after base change its weights on $T_{y}\Orbit_{y}$ occur
in opposite root pairs. The action of $\mathbf{C}_{y}(F)$ therefore
factors through the finite quotient $\mathbf{C}_{y}(F)/\mathbf{C}^{\circ}_{y}(F)$,
so the density $\left|\alpha_{y}\right|_{F}$ is $\mathbf{C}_{y}(F)$-invariant.
Since $\Orbit_{y}\simeq\bH_{s}(F)/\mathbf{C}_{y}(F)$, $\left|\alpha_{y}\right|_{F}$
extends uniquely, by transporting along the conjugation action, to
a $\bH_{s}(F)$-invariant density on $\Orbit_{y}$, inducing a measure
$\mu_{y}$. Concretely, writing $c_{g}(x)=gxg^{-1}$, the density
$\rho_{y}$ defining $\mu_{y}$ is characterized by $(dc_{g}|_{y})^{*}\rho_{y}(gyg^{-1})=\left|\alpha_{y}\right|_{F}$,
for $g\in\bH_{s}(F)$.

\textbf{Step 3}: Define $J_{F}(g,y)\in F^{\times}$ by $(dc_{g}|_{y})^{*}\left(\omega_{\Orbit_{y}}|_{gyg^{-1}}\right)=J_{F}(g,y)\,\omega_{\Orbit_{y}}|_{y}$.
Then $J$ is a definable function on $H\times_{S}H^{\mathrm{rss}}$,
and $\left|J_{F}(g,y)\right|_{F}$ is determined by $(gyg^{-1},y)$,
by Step 2. For $x=gyg^{-1}\in\Orbit_{y}$, set $\beta_{F}(x,y):=\val_{F}J_{F}(g,y)$.
Then $\beta:W\rightarrow\Z$ is definable, and, 
\[
d\mu_{y}(x)=q^{\beta_{F}(x,y)}_{F}\,d\mu^{\mathrm{can}}_{\Orbit_{y}}(x).
\]
A semisimple geometric conjugacy orbit is closed in the algebraic
group. Inside its $F$-points, every rational conjugacy orbit is open:
the orbit map is a submersion, since stabilizers are smooth in characteristic
zero. The complement of any one rational orbit is a union of other
open rational orbits. Thus $\Orbit_{y}$ is also closed in $\bH_{s}(F)$.
Since the restriction of $f_{F}(\cdot,a)$ to $\bH_{s}(F)$ is compactly
supported and locally constant, the integral $\int_{\Orbit_{y}}f_{F}(x,a)\,d\mu_{y}(x)$
converges absolutely for each $(y,a)\in H^{\mathrm{rss}}_{F}\times A_{F}$.
Since $(q^{\beta_{F}}_{F})_{F}\in\mathcal{C}(W)$, it follows from
\cref{prop:canonical-integration} that 
\[
\int_{\Orbit_{y}}f_{F}(x,a)\,d\mu_{y}(x)=\int_{\Orbit_{y}}\left(q^{\beta_{F}(x,y)}_{F}f_{F}(x,a)\right)\,d\mu^{\mathrm{can}}_{\Orbit_{y}}(x),
\]
is constructible on $P\times A$ for each piece $P$ of $H^{\mathrm{rss}}$.
Combining the finitely many pieces proves the lemma. 
\end{proof}

\subsection{Compact elements}

In the embedded family of \cref{thm:families-intro}, fix $F\in\mathrm{Loc}_{0}$
and $s\in\bS(F)$, and write $\bG=\mathscr{G}_{s}$, $G=\bG(F)$,
and $\rho=\rho_{s}$. The definability assertions below are uniform
in $F$ and $s$. An element $x\in G$ is called \emph{compact} if
it generates a relatively compact subgroup. Equivalently, the characteristic
polynomial of $\rho(x)$ has coefficients in $\OO_{F}$ and a unit
constant term. Indeed, compactness forces all eigenvalues to be units.
Conversely, by the Cayley--Hamilton theorem and induction, the set
$\{\rho(x)^{m}:m\in\mathbb{Z}\}$ is bounded in $M_{d}(F)$. In particular,
the set $G_{c}$ of compact elements is closed, invariant under conjugation
and definable. Let $Y=Y_{F,s}$ be the set of regular semisimple compact
elements. It is also definable: $\rho(y)$ is semisimple if and only
if it is annihilated by a square-free polynomial in $F[t]$ of degree
at most $d$, and $y$ is regular if and only if $\dim\ker(\Ad(y)-1)$
equals the absolute rank of $\bG$.

We will use the following standard results, all in characteristic zero. 
\begin{enumerate}[series=representationfacts,label=\textup{(\ensuremath{\mathrm{R}_{\arabic*}})},ref=\ensuremath{\mathrm{R}_{\arabic*}}]
\item \label{item:R-howe} Howe's conjecture (Clozel) \cite[Theorem~1]{MR986793}.
If $\Omega\subset G$ is compact and $J$ is compact
open, the restrictions to $\Hk(G,J)$ of invariant distributions supported
in the closure of $\bigcup_{g\in G}g\Omega g^{-1}$ span a finite-dimensional
space.
\end{enumerate}
Take $\Omega$ to be the union of representatives of the finitely many
conjugacy classes of maximal compact subgroups of $G$. For each compact
element $x$, the compact subgroup $\overline{\langle x\rangle}$ lies in a
maximal compact subgroup, so the conjugacy class of $x$ meets $\Omega$.
Applying \cref{item:R-howe} with this $\Omega$ gives the following consequence.
\begin{enumerate}[label=\textup{(\ensuremath{\mathrm{R}_{1}'})},ref=\ensuremath{\mathrm{R}_{1}'}]
\item \label{item:R-compact-orbits} For every compact open subgroup $J$ of $G$,
the restrictions to $\Hk(G,J)$ of the orbital integrals $O_{y}$,
$y\in Y$, span a finite-dimensional space.
\end{enumerate}
\begin{enumerate}[resume*=representationfacts]
\item \label{item:R-hc} (Harish-Chandra \cite[Theorem~16.3]{MR1702257})
The character $\Theta_{\pi}$ of an irreducible representation $\pi$
is given by a locally integrable function, locally constant on the
regular semisimple set, and Weyl's integration formula applies to
$\Theta_{\pi}(f)$. 
\item \label{item:R-deligne} (Deligne \cite[p.~A156]{MR0425033}) If $\pi$
is supercuspidal then $\Theta_{\pi}$ is supported in the inverse
image of the set of compact elements of $G_{\mathrm{ad}}=\bG_{\mathrm{ad}}(F)$. 
\end{enumerate}
Since $F$ has characteristic zero, the $F$-morphism $\bG\longrightarrow\bG_{\mathrm{ad}}$
is smooth. Hence the image of the induced homomorphism $q:G\to G_{\mathrm{ad}}$
on $F$-points, identified with $G/Z$ with the quotient topology,
is an open, hence closed, subgroup of $G_{\mathrm{ad}}$. Thus the
inverse image in \cref{item:R-deligne} is the set of elements of
$G$ which are compact modulo $Z$. By \cref{item:R-hc}, the function
$\Theta_{\pi}$ vanishes at every regular semisimple element which
is not compact modulo $Z$. 

\begin{lemma}[Supercuspidal characters at a fixed level]
\label{lem:fixed-level-characters}
For every compact open subgroup $J$ of $G$, there exist compact regular
semisimple elements $y_{1},\dots,y_{m}\in G$ such that for every
irreducible supercuspidal representation $\pi$ of $G$ there are scalars
$c_{1}(\pi),\dots,c_{m}(\pi)$ satisfying
\[
\Theta_{\pi}(f)=\sum_{i=1}^{m}c_{i}(\pi)\,O_{y_{i}}(f)
\]
for every $f\in\Hk(G,J)$ supported in $G^{1}$.
\end{lemma}

\begin{proof}
By \cref{item:R-compact-orbits}, choose $y_{1},\dots,y_{m}\in Y$ so that
their orbital integrals span the restrictions of all $O_y$, $y\in Y$,
to the space of such $f$. Suppose that $O_{y_i}(f)=0$ for every $i$.
By Harish-Chandra's regularity theorem \cref{item:R-hc}, Weyl's
integration formula applies to $\Theta_{\pi}(f)$. Only regular
semisimple elements of $G^{1}$ can contribute, since $G^{1}$ is normal
and $f$ is supported there. By Deligne's theorem \cref{item:R-deligne},
such an element must be compact modulo $Z$, hence compact because
$Z\cap G^{1}$ is compact. On each contributing orbit the measure in
Weyl's formula is proportional to the one defining $O_y$, by uniqueness
of invariant measure up to scalar. Thus $\Theta_{\pi}(f)=0$, so
$\Theta_{\pi}$ factors through $f\mapsto (O_{y_i}(f))_{i=1}^{m}$.
The asserted linear combination follows.
\end{proof}

\subsection{Positive shells}

Let $\mathscr{G}$ be a reductive group scheme over a finite-type
$\Z$-scheme $\bS$, as in \cref{thm:families-intro}. On each fiber
$G_{s}=\mathscr{G}_{s}(F)$ of the family, for $s\in S_{F}$, we apply
the argument of \cref{lem:orbital-density} to the action of $G_{s}$
on itself by left translation. Let $\left(\omega_{G_{s}}\right)_{s\in S_{F}}$
be the family of top forms defining the embedded canonical measures.
By translating their value at the identity $e_{s}$, we obtain a family
of left-invariant top forms $\left(\eta_{s}\right)_{s\in S_{F}}$,
with $\eta_{s}|_{g}=(d_{e}L_{g})_{*}(\omega_{G_{s}}|_{e_{s}})$, and
hence Haar measures $dg=\left|\eta_{s}\right|_{F}$ on $G_{s}$. As
in Step 3, the density of $\left|\eta_{s}\right|_{F}$ relative to
$d\mu^{\mathrm{can}}_{G_{s}}$ is $\left|\frac{(d_{e}L_{g})_{*}(\omega_{G_{s}}|_{e_{s}})}{\omega_{G_{s}}|_{g}}\right|_{F}$,
which is an absolute value of a definable function, thus constructible.
Thus, \cref{prop:canonical-integration} applies to integration against
$dg$ on every absolutely integrable fiber.

The following constructions are uniform in $F$ and $s$. When working
on a fixed fiber, we write $G,Z,\rho$ for $G_{s},Z_{s},\rho_{s}$
and abbreviate $\ell_{\rho_{s}}$ and $\ell_{\rho_{s},Z_{s}}$ to
$\ell$ and $\ell_{Z}$. For $r\ge0$, let 
\[
D_{r}=\{g\in G:\ell(g)=\ell_{Z}(g)=r\}.
\]
This set is definable and both open and closed, since $\ell$ and
$\ell_{Z}$ are locally constant. It is compact because it is a closed
subset of $\matball{\rho}(r)$, and every coset $gZ$ with $\ell_{Z}(g)=r$
meets $D_{r}$. Put $v_{n}=\vol(\matcong n{\rho})$. For $g\in G$,
define 
\[
\heckefn ng=1_{\matcong n{\rho}}*\delta_{g}*1_{\matcong n{\rho}}\ \in\ \Hk(G,\matcong n{\rho}).
\]
The volume $v_{n}$ is constructible by \cref{thm:integration} and
\cref{prop:canonical-integration}. Define 
\begin{equation}
\shelltest_{n,r}=\int_{D_{r}}\heckefn ng*\heckefn n{g^{-1}}\,dg\ \in\ \Hk(G,\matcong n{\rho}),\label{eq:shell}
\end{equation}
where $dg$ is the Haar measure fixed above. Explicitly, 
\begin{equation}
\shelltest_{n,r}(x)=v_{n}\int_{D_{r}}\int_{\matcong n{\rho}\times\matcong n{\rho}}1_{\matcong n{\rho}}(g^{-1}u^{-1}xv^{-1}g)\,du\,dv\,dg.\label{eq:shell-explicit}
\end{equation}

\begin{lemma}[The positive shells]
 \label{lem:shell} The functions $\shelltest_{n,r}$ form a constructible
family of bi-$\matcong n{\rho}$-invariant functions, supported in
$\matball{\rho}(2r)\cap G^{1}$. 
\end{lemma}

\begin{proof}
With the unnormalized indicators, \eqref{eq:shell-explicit} is the
product of the constructible volume $v_{n}$ and a constructible integral
by \cref{prop:canonical-integration}. The integrand vanishes unless
$x\in\matcong n{\rho}g\matcong n{\rho}g^{-1}\matcong n{\rho}\subset\matball{\rho}(2r)$.
For every rational character $\chi$ of $\bG$ we have $\abs{\chi}=1$
on $\matcong n{\rho}$, so that $\abs{\chi(x)}=1$ for $x$ in the
support, so $x\in G^{1}$. 
\end{proof}

\section{The support theorem}

\label{sec:char0}

In this section we prove \cref{thm:families-intro} and derive \cref{thm:main},
and \cref{thm:rational-intro} from it.

Let $\bS\hookrightarrow\A^{D}_{\Z}$ be a closed subscheme, let $\mathscr{G}\to\bS$
be a reductive group scheme, and fix a closed embedding $\rho:\mathscr{G}\hookrightarrow\GL_{d,\bS}$
as in \cref{thm:families-intro}. Write $S_{F}=\bS(F)$ and $S=(S_{F})_{F\in\Loc}$.
For $F\in\Loc$ and $s\in\bS(F)$, identify $G_{s}=\mathscr{G}_{s}(F)$
with its image under $\rho_{s}$. Write $Z_{s}$ for its center, $\matcong n{\rho_{s}}$
for its matrix congruence subgroup of level $n$, and $\matball{\rho_{s}}(R)$
for the matrix ball defined by $\ell_{\rho_{s}}$.

First, we state a local rigidity property of reductive groups, which
says that a small perturbation of a reductive group $\mathscr{G}_{s_{0}}(F)$
must be isomorphic to it, and in fact conjugate to it by a matrix
sufficiently close to $e\in\GL_{d}(F)$.

This lets us bound the vanishing threshold for the matrix coefficients
uniformly on each parameter ball. 
\begin{lemma}[Local conjugacy of reductive subgroups]
\label{lem:local-integral-conjugacy}Fix $F\in\Loc$, $s_{0}\in\bS(F)$,
and an open neighborhood $\Omega$ of $e$ in $\GL_{d}(F)$.
There is an open neighborhood $U$ of $s_{0}$ in $\bS(F)$ such that
for every $s\in U$ there exists $a\in\Omega$ satisfying
\[
\rho_{s}(\mathscr{G}_{s})=a\rho_{s_{0}}(\mathscr{G}_{s_{0}})a^{-1}
\]
as closed subgroups of $\GL_{d,F}$. In particular, one may take
$\Omega=\GL_{d}(\OO_{F})$.
\end{lemma}

\begin{proof}
Put $\bS'=\bS\times_{\Z}F$, $\bH_{s}=\rho_{s}(\mathscr{G}_{s})$ and
$\bH=\bH_{s_{0}}$. The group scheme $\mathscr{G}_{\bS'}$ acts through
$\rho$ by left translation on the second factor of
$\mathbf{X}=\bS'\times_{F}(\GL_{d,F}/\bH)$. Since $F$ has characteristic
zero, $\bS'$ is a $\Q$-scheme. The reductive group scheme
$\mathscr{G}_{\bS'}$ acts on the smooth $\bS'$-scheme $\mathbf{X}$, so
its fixed-point scheme is smooth over $\bS'$ by
\cite[Corollary~9.8(4)]{SGA3XII}. Let $\mathbf{T}\subseteq\bS'\times_{F}\GL_{d,F}$
be the inverse image of this fixed-point scheme under the smooth quotient
map $(s,a)\mapsto(s,a\bH)$. Thus $(s,a)\in\mathbf{T}(F)$ precisely when
$a^{-1}\bH_{s}a\subseteq\bH$.

The projection $\mathbf{T}\to\bS'$ is smooth, so it induces an open map
on $F$-points \cite[Proposition~3.5.73(ii)]{MR3729254}. Since
$(s_{0},e)\in\mathbf{T}(F)$, the image of
$\mathbf{T}(F)\cap(\bS(F)\times\Omega)$ under projection to $\bS(F)$
contains an open neighborhood $U$ of $s_{0}$. Thus for each $s\in U$
there is an $a\in\Omega$ with $a^{-1}\bH_{s}a\subseteq\bH$. After
shrinking $U$, the two groups have the same dimension by local constancy
of the fiber dimension of $\mathscr{G}$. Since they are connected,
the inclusion is an equality.
\end{proof}

\begin{proof}[Proof of \cref{thm:families-intro}]
Write $Y_{F,s}$ for the set of regular semisimple compact elements
of $G_{s}$. By \cref{lem:orbital-density,lem:shell}, the functions
$u_{F}(y,s,n,r)=O_{y}(\shelltest_{s,n,r})$ form a constructible family
for $y\in Y_{F,s}$, $n\ge1$, and $r\ge0$. For fixed $F,s,n$, put
\[
W_{F,s,n}=\Span_{\C}\bigl\{(O_{y}(\shelltest_{s,n,r}))_{r\ge0}:y\in Y_{F,s}\bigr\}.
\]
All the test functions $\shelltest_{s,n,r}$ belong to $\Hk(G_{s},\matcong n{\rho_{s}})$,
so \cref{item:R-compact-orbits} shows that $W_{F,s,n}$ is finite-dimensional.

Applying \cref{prop:threshold} with base $Z=S\times\Z_{\ge1}$ and
$z=(s,n)$, we obtain a definable function $R:Z\to\Z_{\geq0}$ such
that every eventually vanishing sequence in $W_{F,s,n}$ vanishes
for $r\ge R_{F}(s,n)$. The same proposition gives an integer $A\ge0$
and a definable function $\gamma:S\to\Z$ such that 
\[
R_{F}(s,n)\le An+\gamma_{F}(s).
\]
Replacing $\gamma$ with $\max\left\{ 0,\gamma\right\} $ we may assume
$\gamma\geq0$.

Fix $F,s,n$. For the next part, we set $\bG=\mathscr{G}_{s}$, $G=G_{s}$,
$\rho=\rho_{s}$, $Z=Z_{s}$, $Y=Y_{F,s}$, $W=W_{F,s,n}$, and $R=R_{F}(s,n)$.

Let $\pi$ be an irreducible supercuspidal representation of $G$.
Unramified twisting preserves the $\matcong n{\rho}$-fixed vectors and
the supports of matrix coefficients. After such a twist, we may assume
that $\pi$ has unitary central character and hence, by supercuspidality,
is unitarizable. Fix a $G$-invariant Hermitian inner product on $\pi$.
Let $\mathcal{E}_{n}=\pi(e_{\matcong n{\rho}})$ be the orthogonal
projection onto $\pi^{\matcong n{\rho}}$; it has finite rank. Set
$a_{\pi}(r)=\Theta_{\pi}(\shelltest_{n,r})$. We have $\pi(\heckefn ng)=v^{2}_{n}\mathcal{E}_{n}\pi(g)\mathcal{E}_{n}$.
Since $\pi(g^{-1})=\pi(g)^{*}$ and $D_{r}$ is compact, 
\begin{equation}
a_{\pi}(r)=v^{4}_{n}\int_{D_{r}}\norm{\mathcal{E}_{n}\pi(g)\mathcal{E}_{n}}^{2}_{\mathrm{HS}}\,dg\ge0.\label{eq:HS}
\end{equation}
The matrix coefficients of $\pi$ with respect to bases of $\pi^{\matcong n{\rho}}$
and $(\pi^{\vee})^{\matcong n{\rho}}$ are compactly supported modulo
$Z$. Hence $a_{\pi}(r)=0$ for $r$ large.

By \cref{lem:shell}, each $\shelltest_{n,r}$ is supported in $G^{1}$.
Applying \cref{lem:fixed-level-characters} with $J=\matcong n{\rho}$,
we see that $a_{\pi}$ is a linear combination of the sequences
$r\mapsto O_{y_i}(\shelltest_{n,r})$, $y_i\in Y$, with coefficients
independent of $r$. Thus $a_{\pi}\in W$. Since $a_{\pi}$ eventually
vanishes, \cref{prop:threshold} gives $a_{\pi}(r)=0$ for all $r\ge R$.

By \eqref{eq:HS}, the integrand being continuous and non-negative
and $D_{r}$ open, $\mathcal{E}_{n}\pi(g)\mathcal{E}_{n}=0$ for $g\in D_{r}$,
$r\ge R$. Since every coset $gZ$ with $\ell_{Z}(g)=r$ meets $D_{r}$
and $Z$ acts by scalars, $\mathcal{E}_{n}\pi(g)\mathcal{E}_{n}=0$
whenever $\ell_{Z}(g)\ge R$. Finally, $\langle\pi(g)v,v^{\vee}\rangle=\langle\mathcal{E}_{n}\pi(g)\mathcal{E}_{n}v,v^{\vee}\rangle=0$
for all $v\in\pi^{\matcong n{\rho}}$, $v^{\vee}\in(\pi^{\vee})^{\matcong n{\rho}}$
and $\ell_{Z}(g)\ge R$.

The threshold obtained above satisfies $R_{F}(s,n)\leq An+\gamma_{F}(s)$.
We now use \cref{lem:local-integral-conjugacy} to make this bound
uniform on each parameter ball $\norm s_{F}\le q^{k}_{F}$.

Write $s\sim_{F}s'$ if $aG_{s}a^{-1}=G_{s'}$ for some $a\in\GL_{d}(\OO_{F})$.
Conjugation by $a$ identifies $\matcong n{\rho_{s}}$ with $\matcong n{\rho_{s'}}$
and preserves the length. It therefore carries $Z_{s}\matball{\rho_{s}}(r)$
onto $Z_{s'}\matball{\rho_{s'}}(r)$ for every $r\ge0$. If $m$ is
a matrix coefficient of a supercuspidal representation of $G_{s}$,
the corresponding matrix coefficient in $G_{s'}$ is $g\mapsto m(a^{-1}ga)$.
Thus any support bound for $G_{s'}$ is also valid for $G_{s}$, at
the same level $n$.

With this discussion in mind, for $F\in\Loc$, we set 
\[
\widetilde{\gamma}_{F}(s)=\min\{\gamma_{F}(s'):s'\in S_{F},\ s'\sim_{F}s\}.
\]
Then $\widetilde{\gamma}:S\rightarrow\Z$ is definable. By the preceding
paragraph, $An+\widetilde{\gamma}_{F}(s)$ is a support bound for
$G_{s}$ for every $n\ge1$. By \cref{lem:local-integral-conjugacy},
$\widetilde{\gamma}_{F}(s)$ is locally constant on $S_{F}$.

For $k\ge0$, let $S_{F,\le k}=\left\{ s\in S_{F}:\norm s_{F}\le q^{k}_{F}\right\} $.
This set is compact, and hence $\widetilde{\gamma}_{F}$ takes only
finitely many values on $S_{F,\le k}$. Set 
\[
\beta_{F}(k)=\max\bigl(\{0\}\cup\{\widetilde{\gamma}_{F}(s):s\in S_{F,\le k}\}\bigr).
\]
Apply \cref{lem:vg-bound} to the definable function $m\mapsto\beta_{F}(m-1)$
on $\Z_{\ge1}$. This gives integers $B,C_{0}\ge0$ and $N\ge1$ such
that $\beta_{F}(m-1)\le Bm+C_{0}+\val_{F}(N)$ for every $F$ and
$m\ge1$. Taking $m=k+1$ and $C=C_{0}+B$ gives, for every $k\ge0$,
\[
\beta_{F}(k)\le Bk+C+\val_{F}(N).
\]
For $s\in S_{F,\le k}$ the matrix support radius is therefore at
most $An+\beta_{F}(k)$. Thus every matrix coefficient $c$ between
$\matcong n{\rho_{s}}$-fixed vectors of an irreducible supercuspidal
representation of $G_{s}$ satisfies
\begin{equation}
\supp(c)\subseteq Z_{s}\matball{\rho_{s}}\bigl(An+Bk+C+\val_{F}(N)\bigr).\label{eq:family-matrix-support}
\end{equation}
For the induced Lie lattice, \eqref{eq:matrix-adjoint-balls} gives
$Z_{s}\matball{\rho_{s}}(R)\subseteq\ball_{s}(2R)$. Hence \cref{thm:families-intro}
holds with constants $2A,2B,2C,N^{2}$. 
\end{proof}

\begin{proof}[Proof of \cref{thm:main}]
Fix $F$ and $\bG$. Before choosing a tower, fix a closed embedding
$\rho:\bG\hookrightarrow\GL_{d}$. The defining equations and group
law descend to a finitely generated $\Z$-subalgebra of $F$. After
localizing this subalgebra, we obtain a reductive group scheme $\mathscr{G}$
over its spectrum $\bS$ and a closed immersion $\mathscr{G}\hookrightarrow\GL_{d,\bS}$.
The inclusion of the subalgebra into $F$ gives a point $s\in\bS(F)$
whose embedded fiber is $(\bG,\rho)$. Choose a closed embedding $\bS\hookrightarrow\A^{D}_{\Z}$
and $k\ge0$ such that $\norm s_{F}\le q^{k}_{F}$. By \eqref{eq:family-matrix-support},
coefficients between $\matcong m{\rho}$-fixed vectors are supported
in $Z\matball{\rho}(Am+c_{0})$, where $c_{0}=Bk+C+\val_{F}(N)$ is
fixed. For the given standard tower, choose $\delta$ as in \eqref{eq:standard-matrix-comparison}.
Vectors fixed by $K_{n}$ in $\pi$ and $\pi^{\vee}$ are fixed by
$\matcong{n+\delta}{\rho}$. Let $e$ be the lattice comparison constant
in \eqref{eq:matrix-adjoint-balls}. Taking $m=n+\delta$ in the matrix
bound and applying this inclusion with $R=A(n+\delta)+c_{0}$ gives
support in $\ball\bigl(2An+2A\delta+2c_{0}+2e\bigr)$. Thus $a=\max\{1,2A\}$
works for every standard tower in $G$, with $b_{K_{\bullet}}=2A\delta+2c_{0}+2e$
for the chosen tower. 
\end{proof}

\begin{corollary}[A fixed model over a number field]
\label{cor:rational-model} Fix a connected reductive group $\bG$ over
a number field $k$ and a closed embedding $\rho:\bG\hookrightarrow\GL_{d,k}$.
There are constants $A,B\ge0$ and an integer $N\ge1$, depending only
on $(k,\bG,\rho)$, with the following property. For every finite place
$v$ of $k$ and every finite extension $F/k_{v}$, put $G=\bG(F)$ and
$Z=Z(G)$, and define $\matcong n{\rho}$ and $\matball{\rho}(R)$ using
the base change of $\rho$ to $F$. For every $n\ge1$, every irreducible
supercuspidal representation $\pi$ of $G$, and all
$w\in\pi^{\matcong n{\rho}}$ and $w^{\vee}\in(\pi^{\vee})^{\matcong n{\rho}}$,
one has
\[
\supp\bigl(g\mapsto\langle\pi(g)w,w^{\vee}\rangle\bigr)
\subset Z\matball{\rho}\bigl(An+B+\val_{F}(N)\bigr).
\]
\end{corollary}

\begin{proof}
Choose an integer $N_{1}\ge1$ such that $(\bG,\rho)$ extends to a
reductive group scheme with a closed embedding in $\GL_{d}$ over
$\bS=\operatorname{Spec}\OO_{k}[1/N_{1}]$. Choose generators
$\alpha_{1},\ldots,\alpha_{t}$ of $\OO_{k}$ as a $\Z$-algebra. Together
with $1/N_{1}$ they give a closed embedding $\bS\hookrightarrow\A^{t+1}_{\Z}$.
For every finite place $v$ and every finite extension $F/k_{v}$, the
corresponding point of $\bS(F)$ has coordinates
$(\alpha_{1},\ldots,\alpha_{t},1/N_{1})$. Since
$|\alpha_{i}|_{F}\le1$ for all $i$, this point has logarithmic height
$\val_{F}(N_{1})$.

Write $A_{0},B_{0},C_{0},N_{0}$ for the constants in
\eqref{eq:family-matrix-support}. At this height, the matrix support
radius modulo the center is at most
\[
A_{0}n+B_{0}\val_{F}(N_{1})+C_{0}+\val_{F}(N_{0})
=A_{0}n+C_{0}+\val_{F}(N_{0}N_{1}^{B_{0}}).
\]
Take $A=A_{0}$, $B=C_{0}$, and $N=N_{0}N_{1}^{B_{0}}$.
\end{proof}

\begin{remark}[Scope of uniformity]
 \label{rem:char0} The field-dependent term in \cref{cor:rational-model}
vanishes outside a finite set of residue characteristics and is controlled
by the absolute ramification index at the remaining ones. The argument
does not remove that term. The proof of the support bound takes place
entirely in characteristic zero. 
\end{remark}

\section{Hecke-algebra bounds and applications}

\label{sec:limit}

We now turn to applications of the support bound. Once the relevant
support bounds are known, the Hecke-algebra estimate below and the
stabilization results of \cref{sec:stab,sec:pairing,sec:adjunction}
follow over non-archimedean local fields of either characteristic.
In characteristic zero, \cref{thm:main} gives the required support
bounds for $G$ and its Levi subgroups, and \eqref{eq:family-matrix-support}
gives the uniform bound in families. For the degree bounds and the
applications over number fields, we use \cite{MR3940435} in its
characteristic-zero setting.

\subsection{The adapted tower}

\label{subsec:adapted-filtrations} We use the filtrations and conventions
of Schneider--Stuhler \cite[Chapter~I, \S2]{MR1471867}. For background
on Bruhat--Tits theory we refer to \cite{MR4520154}. Fix a maximal
$F$-split torus $\adtor$, a special point $x$ in its apartment
in the reduced building, and a minimal parabolic $P_{0}$ containing
$\adtor$. Let $K_{0}=\operatorname{Stab}_{G}(x)\cap G^{1}$. It is
a special maximal compact subgroup. For $n\ge1$, put $K_{n}=U^{(n-1)}_{x}$.
These subgroups form the adapted tower. For a subgroup $H\subseteq G$
write $K_{n,H}=K_{n}\cap H$. For the filtration of the compact subgroup
of $Z_{G}(\adtor)$, fix the one constructed in the proof of \cite[Proposition~I.2.6]{MR1471867}.
We use this same centralizer filtration and the inherited root-group
valuations in every Levi subgroup containing $\adtor$. Thus the induced
filtration on a Levi $L$ is $(K_{n,L})_{n\ge1}$.

Write $\Phi=\Phi(G,\adtor)$ for the relative root system, which need
not be reduced, and $W=W(G,\adtor)$ for the relative Weyl group.
Let $\Phiind$ denote the set of non-divisible roots. Products of
root subgroups are indexed by $\Phiind$; the group attached to a
multipliable root includes its double. Every element of $W$ has a
representative in $N_{G}(\adtor)\cap K_{0}$. Fix such representatives.
They normalize $K_{n}$. The Iwasawa decomposition holds with $K_{0}$
for every parabolic containing $\adtor$.

The groups $K_{n}$, $n\ge1$, are normal in $K_{0}$ and form a neighborhood
basis. For every parabolic $P=MU$ with $M$ containing $\adtor$,
with opposite $\bar{P}=M\bar{U}$, there is Iwahori decomposition
\[
K_{n}=K_{n,U}K_{n,M}K_{n,\bar{U}}=K_{n,\bar{U}}K_{n,M}K_{n,U}.
\]
The same holds inside every Levi subgroup.

Write $H=Z_{G}(\adtor)$ and let $H_{r}$ be its chosen compact-subgroup
filtration. For a real index $r$, the symbol $r+$ means a strict
inequality: $U_{\alpha,r+}=\bigcup_{s>r}U_{\alpha,s}$, and similarly
for $H_{r+}$. Choose an origin for the apartment and put 
\[
d_{\alpha}(n)=n-1-\alpha(x),\qquad R_{\alpha,n}=K_{n,U_{\alpha}}\quad(\alpha\in\Phiind).
\]
The root product formula is 
\begin{equation}
R_{\alpha,n}=\begin{cases}
U_{\alpha,d_{\alpha}(n)+}, & 2\alpha\notin\Phi,\\
U_{\alpha,d_{\alpha}(n)+}\,U_{2\alpha,d_{2\alpha}(n)+}, & 2\alpha\in\Phi,
\end{cases}\qquad K_{n,H}=H_{(n-1)+}.\label{eq:root-depths}
\end{equation}
For any ordering, multiplication of $H_{(n-1)+}$ and the $R_{\alpha,n}$
is a bijective analytic product parametrization of $K_{n}$. Conjugation
satisfies 
\begin{equation}
tU_{\alpha,r}t^{-1}=U_{\alpha,r+\val(\alpha(t))}\qquad(t\in\adtor(F)),\label{eq:root-translation}
\end{equation}
with the same convention for strict indices. The analogous ordered
products for $K_{0}\cap U$, and for the intersection of an Iwahori
subgroup with $U$, have fixed root indices independent of $n$.

For fixed $t\in\adtor(F)$ there is an integer $D(t)\ge0$ such that
\begin{equation}
t^{\pm r}K_{n+D(t)r}t^{\mp r}\subset K_{n}\qquad(n\ge1,\ r\ge0).\label{eq:conj-shift}
\end{equation}

After a fixed shift, the adapted tower is standard in the sense of
\cref{subsec:towers-length}. To see this, choose a chart $\kappa$
over $F$ that is étale at $1$, with $\kappa(1)=0$ and $d\kappa_{1}=\mathrm{id}$,
and a finite splitting field $E$ with valuation extending $\val_{F}$.
Choose a split maximal torus over $E$ whose apartment contains the
image of $x$. In its root and torus coordinates the group $U^{(n-1)}_{x}$
over $E$ is given by $\val_{E}(u_{\alpha})>n-1-\alpha(x)$ on each
root coordinate and $\val_{E}(t_{i}-1)>n-1$ on the torus coordinates
\cite[Proposition~5.2]{MR2905306}. In the resulting normalized chart
these inequalities give $\varpi^{n}\widetilde{\mathcal{L}}_{x}$ at
all sufficiently large depths for a fixed $\OO_{E}$-lattice $\widetilde{\mathcal{L}}_{x}$
in $\Lie(\bG)(E)$. By the Taylor estimate in \cref{subsec:towers-length},
the same lattice describes all sufficiently deep levels in the chart $\kappa$.
Intersecting with $G$ gives
$\kappa(K_{n})=\varpi^{n}(\widetilde{\mathcal{L}}_{x}\cap\Lie(\bG)(F))$
for all sufficiently large $n$ \cite[Lemma~5.3]{MR2905306}.
Rescaling this lattice as in \cref{subsec:towers-length} gives a standard tower
after a fixed shift.

Fix an auxiliary closed embedding $\rho:\bG\hookrightarrow\GL_{d}$.
Comparison with its matrix tower as in \cref{subsec:towers-length}
gives a fixed integer $b\ge0$ such that, for every Levi $L$ containing
$\adtor$ and every $n\ge1$, 
\begin{equation}
\matcong{n+b}{\rho|_{L}}\subset K_{n,L},\qquad K_{n+b,L}\subset\matcong n{\rho|_{L}}.\label{eq:levi-comparison}
\end{equation}
Indeed, compare the two towers in $G$, enlarge $b$ to cover the
finitely many initial levels, and intersect with $L$.

For each such $L$, write $\adball L(R)$ for the adjoint ball defined
by the tower $(K_{n,L})$ via \cref{eq:adjoint-tower-balls}, and
$\adnorm Ll$ for the corresponding operator norm of $\Ad_{L}(l)$.
The same construction in $L$ gives a standard tail of $(K_{n,L})$.
Shifting the tower leaves its eventual conjugation condition unchanged.
In particular, the adjoint balls of the Levi matrix tower are those
of its induced Lie lattice, and \eqref{eq:levi-comparison} compares
them with $\adball L(R)$ by a shift of at most $2b$ in either direction.

All constants may depend on the fixed field, group, and filtration.
Write $K=K_{0}$ and $K_{M}=M\cap K$. Roots and Weyl groups are relative
to $\adtor$.

Since the adapted tower has a standard tail, its lattice description
gives $[K_{1}:K_{n}]\asymp q^{n\dim G}$ for $n\ge1$.

\medskip{}
\noindent\textbf{Cartan double-coset bound.} Over any non-archimedean local
field, the number of $K_{n}$-double cosets in a single Cartan double
coset satisfies 
\begin{equation}
d_{n}(g):=\#(K_{n}\backslash KgK/K_{n})\ll q^{n\dim G}\qquad(g\in G,\ n\ge1),\label{eq:cartan-double-coset-count}
\end{equation}
with an implied constant independent of $g$ and $n$.
\begin{proof}
Write $P_{0}=MU$, with $M=Z_{G}(\adtor)$, and first take $a\in\adtor(F)$
in the chamber defined by $|\alpha(a)|\le1$ for the roots in $U$.
By the Iwahori factorization $K_{n}=K_{n,U}K_{n,M}K_{n,\bar{U}}$
and the root-group filtrations above, 
\[
aK_{n,U}a^{-1}\subseteq K_{n,U},\qquad aK_{n,M}a^{-1}=K_{n,M},\qquad aK_{n,\bar{U}}a^{-1}\supseteq K_{n,\bar{U}}.
\]
Uniqueness of the factors gives $K_{n}\cap aK_{n}a^{-1}=(aK_{n,U}a^{-1})K_{n,M}K_{n,\bar{U}}$.
Consequently 
\[
\frac{\vol(K_{n}aK_{n})}{\vol(K_{n})}=[K_{n,U}:aK_{n,U}a^{-1}]=\bigl|\det(\Ad(a)|_{\Lie U})\bigr|^{-1},
\]
independently of $n$. Since $K_{n}\triangleleft K$, all $K_{n}$-double
cosets in $KaK$ have the same volume. Thus $d_{n}(a)=d_{1}(a)[K_{1}:K_{n}]$,
and $d_{1}(a)\le[K:K_{1}]^{2}$.

The Cartan decomposition parametrizes $K\backslash G/K$ by the Weyl
orbits in $M/M^{1}$. The torus $\adtor$ is the maximal split central
torus of $M$, so $\adtor(F)$ has finite-index image in $M/M^{1}$
by \cref{subsec:towers-length}. Thus every $K$-double coset has
a representative $at$ with $a$ as above and $t$ in a fixed finite
subset of $M$, using Weyl representatives in $K$ to put $a$ in
the chosen chamber. For $H=K$ or $K_{n}$, the groups $H$ and $tHt^{-1}$
meet in bounded index in both groups, uniformly in $n$ by the lattice
description of the tower. Hence replacing $a$ by $at$ changes $[H:H\cap aHa^{-1}]$
by at most a bounded factor. Since $\vol(HgH)=\vol(H)[H:H\cap gHg^{-1}]$,
taking the ratio for $K$ and $K_{n}$ gives $d_{n}(at)\ll d_{n}(a)$.
This proves \eqref{eq:cartan-double-coset-count}. 
\end{proof}

\subsection{The supercuspidal part of the Hecke algebra}

Put $A_{G}=\mathbf{A}_{G}(F)$ for the split central torus of \cref{subsec:towers-length}.
Fix a discrete cocompact subgroup $\centrallattice\subset A_{G}$
and put $\Gbar=G/\centrallattice$. In this subsection let $K_{\bullet}=(K_{n})_{n\ge1}$
be any standard tower, with adjoint balls $\ball(R)$ as in \cref{subsec:towers-length}.
Write $\overline{K_{n}}$ for the image of $K_{n}$ in $\Gbar$. Assume
that the supercuspidal matrix coefficients between $K_{n}$-fixed
vectors are supported in $\ball(R_{n})$, where $R_{n}=a_{G}n+b_{K_{\bullet}}$
as in \cref{thm:main}. Put 
\[
\IrrscZ G=\{\pi\in\Irrsc G:\omega_{\pi}|_{\centrallattice}=1\},
\]
where $\omega_{\pi}$ is the central character of $\pi$. Let $\Hsc(\Gbar,\overline{K_{n}})$
be the cuspidal part of $\Hk(\Gbar,\overline{K_{n}})$. For $\pi\in\IrrscZ G$,
let $\mu_{\pi}$ be its formal-degree measure on $\Gbar$. Identify
$\pi^{K_{n}}\otimes(\pi^{\vee})^{K_{n}}$ with $\operatorname{End}_{\C}(\pi^{K_{n}})$
by sending $v\otimes\lambda$ to $w\mapsto\lambda(w)v$. Schur orthogonality
gives the algebra isomorphism 
\[
\bigoplus_{\pi\in\IrrscZ G}\operatorname{End}_{\C}(\pi^{K_{n}})\xrightarrow{\ \sim\ }\Hsc(\Gbar,\overline{K_{n}}),\qquad v\otimes\lambda\longmapsto\lambda(\pi(g^{-1})v)\,\mu_{\pi}(dg).
\]

\begin{corollary}[An upper bound for the supercuspidal part of the Hecke algebra]
 \label{cor:hecke-upper} Put $r=\operatorname{rank}_{F}(G/Z)$.
There is a constant $C>0$, depending only on the field, group, standard
tower, and $\centrallattice$, such that 
\[
\begin{aligned}\dim\Hsc(\Gbar,\overline{K_{n}}) & =\sum_{\pi\in\IrrscZ G}(\dim\pi^{K_{n}})^{2}\\
 & \le\#\bigl\{\centrallattice K_{n}gK_{n}:g\in\ball(R_{n})\bigr\}\le C(1+n)^{r}q^{n\dim G}\qquad(n\ge1).
\end{aligned}
\]
\end{corollary}

\begin{proof}
Each density $g\mapsto\lambda(\pi(g^{-1})v)$ is a matrix coefficient
of $\pi^{\vee}$, so the corresponding Hecke measure is supported
on the image of $\ball(R_{n})$ in $\Gbar$. On each $\overline{K_{n}}$-double
coset, biinvariant measures form a one-dimensional space. This gives
the first inequality.

By the tower comparison in \cref{subsec:towers-length}, it suffices
to count double cosets for the adapted tower at level $m=n+d$ in its
ball of radius $R=R_{n}+2d$, for a fixed $d\ge0$.
If $KatK$ meets this ball,
then $\norm{\Ad(a)}\le q^{R+O(1)}$, since $K$ is compact and $t$
ranges over a fixed finite set. Evaluating $\Ad(a)$ on the root spaces
for $\alpha$ and $-\alpha$ gives $\abs{\val(\alpha(a))}\le R+O(1)$
for each relative root. The roots span the dual of the apartment modulo
the split center, and $\centrallattice$ is cocompact in $A_{G}$.
Thus modulo $\centrallattice$ there are $O((1+R)^{r})$ Cartan parameters
meeting the ball. By \eqref{eq:cartan-double-coset-count} each contributes
at most $O(q^{m\dim G})$ double cosets at level $m$. Since $m=n+d$
and $R=O(n)$, this gives the required bound. 
\end{proof}

When $Z$ is compact, \cite[Theorem~6.14]{LS26} applied with compact
open subgroup $K_{1}$ gives a lower bound for $\dim\Hsc(G,K_{n})$.
For general $Z$, the same argument applies to $\Hsc(\Gbar,\overline{K_{n}})$.
Together with \cref{cor:hecke-upper}, this yields 
\[
q^{n\dim G}\ll\dim\Hsc(\Gbar,\overline{K_{n}})\ll(1+n)^{r}q^{n\dim G}
\]
for all sufficiently large $n$, with constants allowed to depend
on the chosen tower and $\centrallattice$. Thus the two bounds differ by a factor $C(1+n)^{r}$.

\subsection{Degree bounds and limit multiplicities}

In this subsection $F$ has characteristic zero, as in \cite{MR3940435}.
We relate \cref{thm:main} to property (PSC) of that paper. For the
fixed embedding $\rho:\bG\hookrightarrow\GL_{d}$ and a connected
reductive subgroup $\bH$ of $\bG$, write $H=\bH(F)$ and put $\matcong n{\rho|_{\bH}}=\matcong n{\rho}\cap H$.
Choose a lattice in $\Lie(\bH)(F)$ and use its operator norm to define
\[
\adball H(R)=\{h\in H:\norm{\Ad_{H}(h)}\le q^{R}\}.
\]
When the chosen lattice is the Lie lattice of a standard tower on
$H$, these are its adjoint balls as in \cref{subsec:towers-length}.
By \cite[Definition~3]{MR3940435}, $H$ has property (PSC) if there
is $c>0$ such that for every quasi-cuspidal representation $\pi$
of $H$ the matrix coefficients between $\matcong n{\rho|_{\bH}}\cap\bH^{\mathrm{der}}(F)$-fixed
vectors are supported in $\adball H(cn)$, $n\ge1$. As noted immediately
after that definition, it suffices to consider irreducible supercuspidal
representations, and (PSC) for $H$ is equivalent to (PSC) for $\bH^{\mathrm{der}}(F)$.
\begin{lemma}[Property (PSC)]
\label{lem:psc} For every connected reductive group $\bH$ over
$F$, the group $\bH(F)$ has property (PSC). 
\end{lemma}

\begin{proof}
By the reduction to the derived group we may assume that $\bH$ is semisimple.
For some fixed $e\ge0$, $(\matcong{j+e}{\rho|_{\bH}})_{j\ge1}$ is a
standard tower by \cref{subsec:towers-length}. Since $\matcong{n+e}{\rho|_{\bH}}
\subseteq\matcong n{\rho|_{\bH}}$, \cref{thm:main} bounds the support
radius by $a_{H}n+O(1)$. Changing the Lie lattice adds a bounded constant,
which can be absorbed in $cn$ for $n\ge1$.
\end{proof}

\begin{remark}[Uniformity over all places]
\label{rem:all-places} Let $\bH$ be a semisimple group over a number
field $k$, with the embedding and lattices defined globally. Then the
constant in \cref{lem:psc} for $\bH(k_{v})$ can be taken independent
of the finite place $v$. Indeed, apply \cref{cor:rational-model} with
$F=k_{v}$. Since $\sup_{v\nmid\infty}\val_{k_{v}}(N)<\infty$, the
matrix support radius modulo the center is at most $cn$ for all $v$
and $n\ge1$. By \eqref{eq:matrix-adjoint-balls}, the adjoint radius for
the lattice induced by $d\rho$ is at most $2cn$. The prescribed global
Lie lattice agrees with the induced one at all but finitely many places.
At the remaining places, the bounded additive changes in radius can be
absorbed by increasing $c$. This complements \cite{MR3940435}, where
(PSC) with a uniform constant is proved outside a finite set of places.
\end{remark}

By \cref{lem:psc} and \cite[Theorem~1]{MR3940435}, every connected reductive
group over $F$ satisfies (BD). Over a number field, the same argument
gives a uniform bound. Outside finitely many places, the semisimple normal
factors of the standard local Levi subgroups occur at integral points of a
fixed finite-type family over the scheme of Borel pairs of a reductive
integral model. Applying \eqref{eq:family-matrix-support} at height zero
on a finite affine cover and then \eqref{eq:matrix-adjoint-balls} gives a
uniform (PSC) constant for these factors. The remaining constants in the
passage to (BD) are uniform outside finitely many places by the proof of
\cite[Theorem~4]{MR3940435}. Taking the maximum with the local constants
at the exceptional places gives uniform (BD) at all finite places.

\begin{corollary}[Limit multiplicities]
\label{cor:limit} Let $\bH$ be a reductive group over a number
field $k$ which satisfies property (TWN) \cite{MR3352530}. Let $S$
be a finite set of places containing all archimedean places, and let
$\Adeles^{S}$ denote the adeles away from $S$. Then the limit multiplicity
property holds for every non-degenerate family of open subgroups of
a fixed open compact subgroup $K^{S}_{0}$ of $\bH(\Adeles^{S})$,
in the sense of \cite{MR3940435}. This applies in particular to quasi-split
classical groups. 
\end{corollary}

\begin{proof}
Apply \cite[Theorem~5]{MR3940435} with $S_{0}=\emptyset$, using
(BD) as established above and the assumed property (TWN). Quasi-split
classical groups satisfy (TWN) by \cite{MR3711830}, as used in \cite{MR3940435}. 
\end{proof}

\section{Bruhat strata}

\label{sec:strata}

We now allow $F$ to have either characteristic. We use the adapted
tower of \cref{subsec:adapted-filtrations} and the Bruhat-filtration
method of \cite[\S3]{MR3940435}. For the standard stabilization,
duality, and adjointness constructions we follow Bushnell \cite{MR1810135};
we retain the estimates that make their thresholds linear in the level.

\subsection{Contracting elements}

For every standard parabolic subgroup $P=MU$ fix an $F$-defined
cocharacter $\lambda_{P}$ of $\adtor$ such that 
\[
\langle\alpha,\lambda_{P}\rangle=0\quad(\alpha\in\Phi_{M}),\qquad\langle\alpha,\lambda_{P}\rangle>0\quad(\alpha\in\Phi_{U}),
\]
and set $t_{P}=\lambda_{P}(\varpi)\in Z(M)\cap\adtor(F)$. Thus $\val(\alpha(t_{P}))=\langle\alpha,\lambda_{P}\rangle\ge1$
for $\alpha\in\Phi_{U}$, i.e., $t_{P}$ contracts $U$. Let $a_{P}=t^{-1}_{P}$
and 
\[
U[r]=a^{r}_{P}(U\cap K)a^{-r}_{P},\qquad r\in\Z_{\ge0}.
\]
This is an increasing sequence of compact open subgroups of $U$,
normalized by $K_{M}=M\cap K$, whose union is $U$. Choose an integer
$D_{P}$ for which \eqref{eq:conj-shift} holds with $t=t_{P}$. By
\eqref{eq:conj-shift}, applied twice, and the normality of $K_{n}$
in $K$, we have 
\begin{equation}
u^{-1}K_{n+2D_{P}r}u\subset K_{n}\qquad\text{for \ensuremath{u\in U[r]} and \ensuremath{n\ge1}.}\label{eq:unipotent-shift}
\end{equation}
Indeed, write $u=a^{r}_{P}u_{0}a^{-r}_{P}$ with $u_{0}\in K$; then
\[
u^{-1}K_{n+2D_{P}r}u\subset a^{r}_{P}u^{-1}_{0}K_{n+D_{P}r}u_{0}a^{-r}_{P}=a^{r}_{P}K_{n+D_{P}r}a^{-r}_{P}\subset K_{n}.
\]

\subsection{Support on a Bruhat stratum}

Let $Q=LN_{Q}$ be another standard parabolic subgroup. Let $w_{1},\dots,w_{s}\in N_{G}(\adtor)\cap K$
be representatives for the double cosets $Q\backslash G/P$, numbered
so that 
\[
Z_{i}=\bigcup_{j>i}Qw_{j}P
\]
is closed for $i=0,\dots,s$. (Any numbering refining the closure
order will do.) In particular, the closure of $Qw_{i}P$ is contained
in $Z_{i-1}=Qw_{i}P\cup Z_{i}$.
\begin{lemma}[Support on a Bruhat stratum]
\label{lem:stratum} There are integers $a_{0}\ge1$, $b_{0}\ge0$
such that for every $i$, $n\ge1$ and $r\ge a_{0}n+b_{0}$ we have
\[
Qw_{i}P\subset Z_{i}K_{n}\cup Qw_{i}U[r]K_{M}.
\]
\end{lemma}

\begin{proof}
Choose a chamber adjacent to $x$ in the direction of the positive
roots of $P_{0}$, and let $\mathcal{I}^{\circ}$ be its connected
Iwahori subgroup. Write $H_{0}=Z_{G}(\adtor)\cap K$ and put $\mathcal{I}=H_{0}\mathcal{I}^{\circ}$.
The compact centralizer normalizes the Iwahori root subgroups. The
root product description at the chamber gives, for a standard $P=MU$,
\begin{equation}
\mathcal{P}:=(P\cap K)\mathcal{I}=(P\cap K)(\bar{U}\cap\mathcal{I})=(\bar{U}\cap\mathcal{I})(P\cap K).\label{eq:parahoric-cell}
\end{equation}
These are the parahoric product decompositions of \cite{MR0327923},
with the root-depth conventions of \cref{subsec:adapted-filtrations}.
The chamber form of Iwasawa decomposition gives $G=P_{0}N_{G}(\adtor)\mathcal{I}$.
Decomposing the relative Weyl group into $W_{L}\backslash W/W_{M}$
therefore gives 
\[
G=\bigcup_{i}Qw_{i}\mathcal{P}.
\]
This is the covering used in \cite[Lemma~2]{MR3940435}; we record
its application to the two arbitrary parabolics $Q$ and $P$.

Suppose that $Qw_{i}P$ meets $Qw_{j}\mathcal{P}$. Then $w_{j}\bar{u}\in Qw_{i}P$
for some $\bar{u}\in\bar{U}\cap\mathcal{I}$. Since $a_{P}$ contracts
$\bar{U}$, $w_{j}a^{r}_{P}\bar{u}a^{-r}_{P}=(w_{j}a^{r}_{P}w^{-1}_{j})w_{j}\bar{u}a^{-r}_{P}\in Qw_{i}P$
tends to $w_{j}$ as $r\to\infty$. Hence $w_{j}$ lies in the closure
of $Qw_{i}P$, and therefore $j\ge i$.

Next we show that $Qw_{i}P\cap Qw_{i}\mathcal{P}=Qw_{i}(P\cap K)$.
Let $R=w^{-1}_{i}Qw_{i}$; it is a parabolic subgroup containing $\adtor$.
We claim that 
\begin{equation}
R\cap\bar{U}MU=(R\cap\bar{U})(R\cap M)(R\cap U).\label{eq:R-cell}
\end{equation}
Work first with the underlying algebraic groups over an algebraic
closure of $F$. Under the big-cell isomorphism $\bar{U}\times M\times U\to\bar{U}MU$,
the product on the right of \eqref{eq:R-cell} is closed in the big
cell. By the root-group decomposition it contains a Zariski-open neighborhood
of the identity in $R$. Since $R$ is irreducible, the product equals
$R\cap\bar{U}MU$. The big-cell isomorphism is defined over $F$,
so uniqueness of the factors gives \eqref{eq:R-cell} on $F$-points.
In particular, $\bar{U}\cap RP=\bar{U}\cap R$. Applied to $w_{i}\bar{u}\in Qw_{i}P$,
this proves the claimed intersection.

We conclude that $Qw_{i}P\subset Z_{i}(\bar{U}\cap\mathcal{I})\cup Qw_{i}(P\cap K)$.
Multiply on the right by $a^{-r}_{P}$. Since $Z_{i}$ and $Qw_{i}P$
are right $P$-invariant and $w_{i}a^{r}_{P}w^{-1}_{i}\in Q$, we
get 
\[
Qw_{i}P\subset Z_{i}\,a^{r}_{P}(\bar{U}\cap\mathcal{I})a^{-r}_{P}\cup Qw_{i}U[r]K_{M}.
\]
By \eqref{eq:root-translation}, conjugation by $a_{P}^{r}$ raises
each root depth of $\bar{U}\cap\mathcal{I}$ by at least $r$, including
the depths for doubled roots. Since the initial depths are fixed,
they exceed the corresponding $d_{\alpha}(n)=n-1-\alpha(x)$ after
conjugation whenever $r\ge n+b_{0}$ for a fixed $b_{0}$. Thus
\eqref{eq:root-depths} gives $a_{P}^{r}(\bar{U}\cap\mathcal{I})a_{P}^{-r}\subset K_{n}$.
\end{proof}

\subsection{Intersections}

% The unipotent radical in the Levi L on the stratum indexed by w and its opposite.
\newcommand{\Uw}{U_{w}}
\newcommand{\oppUw}{\bar{U}_{w}}

We keep the notation of the previous subsection and write $w=w_{i}$.
Let 
\[
H=w^{-1}Qw\cap P,\qquad H_{U}=w^{-1}Qw\cap U,\qquad H_{M}=w^{-1}Qw\cap M,
\]
\[
\Uw=L\cap wUw^{-1},\qquad L'=L\cap wMw^{-1}.
\]

\begin{lemma}[Projection to the Levi subgroup]
\label{lem:intersections}
\begin{enumerate}[label=\textup{(\arabic*)},ref=\arabic*]
\item \label{item:int-H} $H=H_{M}\ltimes H_{U}$, and $H_{M}$ is a parabolic
subgroup of $M$ containing $\adtor$. Moreover $M=H_{M}K_{M}$, and
consequently $QwP=QwUK_{M}$. 
\item \label{item:int-L} $L'\Uw=L\cap wPw^{-1}$ is a parabolic subgroup
of $L$ with Levi part $L'$ and unipotent radical $\Uw$. 
\item \label{item:int-proj} Let $\operatorname{pr}:Q\to L$ be the projection.
For $r\ge0$, $\operatorname{pr}$ maps $w(H_{U}\cap U[r])w^{-1}$
homomorphically onto 
\[
\Uw[r]=b^{-1}_{r}(\Uw\cap K)b_{r},\qquad b_{r}=wt^{r}_{P}w^{-1}.
\]
In particular, it takes the normalized Haar measure of $w(H_{U}\cap U[r])w^{-1}$
to that of $\Uw[r]$. 
\item \label{item:int-b} The element $b_{r}$ centralizes $L'$, and $\abs{\beta(b_{r})}\le q^{-r}$
for every root $\beta$ of $\adtor$ in $\Uw$. In particular, if $\Uw$ is nontrivial,
then $\adnorm L{b_{r}}\ge q^{r}$.
\end{enumerate}
\end{lemma}

\begin{proof}
The root-group decomposition of $w^{-1}Qw\cap P$, with the roots
in $M$ placed before those in $U$, gives $H=H_{M}H_{U}$. The roots
of $H_{M}$ in $M$ are those belonging to $w^{-1}Qw$; they define
a parabolic subgroup of $M$. The Iwasawa decomposition in $M$ gives
$M=H_{M}K_{M}$ and hence $QwP=QwUK_{M}$. Applying the same argument
in $L$ proves \cref{item:int-L}.

The root product decomposition, ordered with the $L$-factors first,
gives 
\[
w(H_{U}\cap K)w^{-1}=(\Uw\cap K)(N_{Q}\cap wUw^{-1}\cap K).
\]
Levi projection is projection onto the first factor. Since $t_{P}$
normalizes $H_{U}$, conjugating this identity by $b^{-1}_{r}$, where
$b_{r}=wt^{r}_{P}w^{-1}$, gives the image $\Uw[r]$ in \cref{item:int-proj}.
The measure assertion follows from surjectivity of this homomorphism
of compact groups.

Finally $b_{r}$ centralizes $L'$ since $t_{P}$ centralizes $M$.
For each relative root $\beta$ in $\Uw$, $\val(\beta(b_{r}))=r\langle w^{-1}\beta,\lambda_{P}\rangle\ge r$.
If $\Uw$ is nontrivial, choose such a $\beta$. On the opposite root space in
$\Lie L$, $\Ad_{L}(b_{r})$ acts by $\beta(b_{r})^{-1}$, giving the
asserted norm bound.
\end{proof}

\section{Averaging on induced representations}

\label{sec:averaging}

From this section onward, assume the \levsupport{} for every standard Levi
subgroup $L$, including $G$. This means that there is a constant
$c_{L}$ such that for every irreducible supercuspidal representation
$\sigma$ of $L$ and $r\ge1$, 
\begin{equation}
\langle\sigma(l)v,v^{\vee}\rangle=0\quad\text{if }l\notin\adball L(c_{L}r),\qquad v\in\sigma^{K_{r,L}},\quad v^{\vee}\in(\sigma^{\vee})^{K_{r,L}}.\label{eq:levi-support}
\end{equation}
In characteristic zero, choose a fixed $d\ge0$ such that $(K_{j+d,L})_{j\ge1}$
is a standard tower. The vectors in \eqref{eq:levi-support} are fixed
by $K_{r+d,L}$, and the shifted tower defines the same adjoint balls.
Applying \cref{thm:main} to $L$ at level $r$ gives an affine radius
in $r$, hence a radius $c_{L}r$ for $r\ge1$, as required in \eqref{eq:levi-support}.
The argument below uses only the \levsupport{}.

Let $Q=LN_{Q}$ and $P=MU$ be standard parabolic subgroups, $\sigma$
an irreducible supercuspidal representation of $L$, and $I=i^{G}_{Q}\sigma$
the normalized induced representation, realized on functions $f:G\to\sigma$
with $f(xg)=\Delta_{Q}(x)^{1/2}\sigma(x)f(g)$ for $x\in Q$, where
$\sigma$ is extended trivially to $N_{Q}$ and $\Delta_{Q}$ is the
modulus function of $Q$. With the notation of \cref{sec:strata},
let 
\[
I_{i}=\{f\in I:f|_{Z_{i}}\equiv0\},\qquad i=0,\dots,s.
\]
Thus $0=I_{0}\subset I_{1}\subset\dots\subset I_{s}=I$ are $P$-invariant
subspaces. For a compact group $C$ we write $e_{C}$ for its normalized
Haar measure, acting on representations.

We use the following averaging identity in $G$ and in its Levi subgroups.
Let $J$ be a compact open subgroup with Iwahori decomposition relative
to $P=MU$. Write $J_{U}=J\cap U$ and similarly for $M$ and $\bar{U}$.
If $t\in Z(M)$ satisfies $t^{-1}J_{\bar{U}}t\subseteq J_{\bar{U}}$, then
for every $J$-fixed vector $v$ in a smooth representation,
\begin{equation}
e_{J}tv=e_{J_{U}}tv=t e_{t^{-1}J_{U}t}v.\label{eq:positive-averaging}
\end{equation}
Indeed, $tv$ is fixed by $J_{M}$ and $J_{\bar{U}}$, and
$e_{J}=e_{J_{U}}e_{J_{M}}e_{J_{\bar{U}}}$.

The next lemma follows the averaging argument of \cite[Lemmas~3--5]{MR3940435}
for the two parabolics $Q$ and $P$.
\begin{lemma}[Averaging over the stratum]
\label{lem:averaging-step} There exists $A\ge1$, depending on the
fixed field, group, filtration, and the choices of $\lambda_{P}$,
such that the following holds for all $Q$, $P$, $\sigma$ and $i$.
Let $r\ge1$ and $f\in I^{K_{r}}_{i}$, and assume that the image
of $f$ in the Jacquet module $I_{U}$ vanishes. Then $e_{U[m]}f\in I_{i-1}$
for all $m\ge Ar$. 
\end{lemma}

\begin{proof}
Write $w=w_{i}$, and keep the notation of \cref{lem:intersections}.
Put $R=a_{0}r+b_{0}$, increasing $a_{0},b_{0}$ to be integers, and
let $m\ge R$. Since $f$ vanishes on $Z_{i}K_{r}$, \cref{lem:stratum}
gives 
\begin{equation}
\supp f\cap QwP\subset QwU[R]K_{M}.\label{eq:support-on-stratum}
\end{equation}
The function $e_{U[m]}f$ lies in $I_{i}$, so it suffices to show
that it vanishes on $QwP=QwUK_{M}$ (\cref{lem:intersections}). We
first show that for every $k\in K_{M}$ 
\begin{equation}
\supp\bigl(u\mapsto f(wuk)\bigr)\subset H_{U}U[R].\label{eq:support-in-U}
\end{equation}
Indeed, suppose that $wuk=xwu_{0}k_{0}$ with $x\in Q$, $u_{0}\in U[R]$
and $k_{0}\in K_{M}$. Then $w^{-1}xw=uk_{1}u^{-1}_{0}k^{-1}_{1}\cdot k_{1}\in H$,
where $k_{1}=kk^{-1}_{0}$. Since $H=H_{M}\ltimes H_{U}$, we get
$k_{1}\in H_{M}$ and $uk_{1}u^{-1}_{0}k^{-1}_{1}\in H_{U}$. As $k_{1}$
normalizes $U[R]$, it follows that $u\in H_{U}U[R]$.

Let $u\in U$ and $k\in K_{M}$. By \eqref{eq:support-in-U}, $(e_{U[m]}f)(wuk)=0$
unless $u\in H_{U}U[m]$, and if $u=hu_{1}$ with $h\in H_{U}$, $u_{1}\in U[m]$,
then $(e_{U[m]}f)(wuk)=\sigma(\operatorname{pr}(whw^{-1}))(e_{U[m]}f)(wk)$,
since $\Delta_{Q}$ is trivial on the unipotent group $wH_{U}w^{-1}$.
Thus it suffices to show that 
\begin{equation}
\int_{U[m]}f(wuk)\,du=0\qquad\text{for all }k\in K_{M}.\label{eq:vanishing-integral}
\end{equation}

Suppose first that $\Uw=L\cap wUw^{-1}$ is nontrivial. Set $r'=r+2D_{P}R$. For $u_{0}\in U[R]$,
$k\in K_{M}$ and $l\in K_{r',L}$ we have $(wu_{0}k)^{-1}l(wu_{0}k)\in k^{-1}u^{-1}_{0}K_{r'}u_{0}k\subset K_{r}$
by \eqref{eq:unipotent-shift}, so that $f(wu_{0}k)\in\sigma^{K_{r',L}}$.
We claim that 
\begin{equation}
e_{\Uw[m]}\sigma^{K_{r',L}}=0\qquad\text{if }m>c_{L}r'.\label{eq:levi-averaging}
\end{equation}
Indeed, let $v\in\sigma^{K_{r',L}}$. By \cref{lem:intersections}\cref{item:int-b}
and \eqref{eq:levi-support}, $\sigma(b_{m})v$ pairs to zero with
$(\sigma^{\vee})^{K_{r',L}}$, so $e_{K_{r',L}}\sigma(b_{m})v=0$.
Apply \eqref{eq:positive-averaging} in $L$ with $J=K_{r',L}$,
parabolic $L'\Uw$, and $t=b_{m}$. By \eqref{eq:root-translation},
$b_{m}^{-1}K_{r',\oppUw}b_{m}\subset K_{r',\oppUw}$ for the opposite
unipotent radical $\oppUw$. Thus $e_{b_{m}^{-1}K_{r',\Uw}b_{m}}v=0$.
Since $b_{m}^{-1}K_{r',\Uw}b_{m}\subset \Uw[m]$, \eqref{eq:levi-averaging} follows.

Now decompose the integral in \eqref{eq:vanishing-integral} according
to $(H_{U}\cap U[m])\backslash U[m]$. By \eqref{eq:support-in-U},
only cosets which meet $U[R]$ contribute (note that $U[m]\cap H_{U}U[R]=(H_{U}\cap U[m])U[R]$).
For $u_{0}\in U[R]$, by \cref{lem:intersections} \cref{item:int-proj},
\[
\int_{H_{U}\cap U[m]}f(whu_{0}k)\,dh=\vol(H_{U}\cap U[m])\,e_{\Uw[m]}f(wu_{0}k)=0
\]
by \eqref{eq:levi-averaging}, provided that $m>c_{L}(r+2D_{P}R)$.

Suppose now that $L\cap wUw^{-1}=\{1\}$. For $k\in K_{M}$ put $\phi(u)=f(wuk)$.
Since $wH_{U}w^{-1}\subset N_{Q}$, $\phi$ is left $H_{U}$-invariant.
The Jacquet image of $f$ is zero, so $e_{U[d]}f=0$ for some $d\ge m$.
Since $k$ normalizes $U[d]$, evaluation at $wk$ gives $\int_{U[d]}\phi(u)\,du=0$.
By \eqref{eq:support-in-U} and $m\ge R$, the support in $U[d]$ lies in
$(H_{U}\cap U[d])U[m]$. This is a finite union of left $U[m]$-cosets
with representatives in $H_{U}$, each contributing $\int_{U[m]}\phi(u)\,du$.
This proves \eqref{eq:vanishing-integral}.

Since $R=a_{0}r+b_{0}\le(a_{0}+b_{0})r$ for $r\ge1$, an integer
$A$ larger than all the finitely many coefficients just obtained
proves the assertion. 
\end{proof}

\begin{proposition}[Averaging on induced representations]
\label{prop:induced} There exists $C$, depending on the fixed field,
group, filtration, and the $\lambda_{P}$, such that for all $Q$,
$P$, $\sigma$ as above and all $n\ge1$, 
\[
e_{U[m]}\bigl(I^{K_{n}}\cap\ker(I\to I_{U})\bigr)=0\qquad\text{for }m\ge Cn.
\]
\end{proposition}

\begin{proof}
Let $f\in I^{K_{n}}$ have zero image in $I_{U}$. Define $m_{0}=0$
and $m_{j}=A(n+2D_{P}m_{j-1})$ for $j=1,\dots,s$, and put $f_{0}=f$
and $f_{j}=e_{U[m_{j}]}f_{j-1}$ for $j\ge1$. Since the groups $U[m]$
increase, $f_{j}=e_{U[m_{j}]}f$ for $j\ge1$. By \eqref{eq:unipotent-shift},
$f_{j-1}$ is fixed by $K_{n+2D_{P}m_{j-1}}$, and its image in $I_{U}$
is zero since averaging over a subgroup of $U$ acts trivially on
$I_{U}$. By induction and \cref{lem:averaging-step}, $f_{j}\in I_{s-j}$,
so that $f_{s}=0$. Clearly $m_{s}\le C'n$ where $C'$ depends only
on $A$, $D_{P}$ and $s\le\abs W$. The same holds for $m\ge m_{s}$
since $e_{U[m]}=e_{U[m]}e_{U[m_{s}]}$. 
\end{proof}

No bound on $\dim I^{K_{n}}$ is used.

\section{From induced representations to all representations}

\label{sec:smooth}

We use qualitative stabilization in the form of \cite[Theorem~1 and \S3.5]{MR1810135},
and then transfer our bound from induced representations without increasing
its exponent. Let $P=MU$ with $M$ containing $\adtor$, let $\bar{P}=M\bar{U}$
be opposite to it, and fix an integer $n\ge1$. The congruence subgroup
$K_{n}$ is totally decomposed relative to the minimal Levi $Z_{G}(\adtor)$,
as required in \cite{MR1810135}, and 
\[
K_{n}=K_{n,U}K_{n,M}K_{n,\bar{U}}.
\]
Let $t\in Z(M)$ satisfy $tK_{n,U}t^{-1}\subset K_{n,U}$, $t^{-1}K_{n,\bar{U}}t\subset K_{n,\bar{U}}$
and $U=\bigcup_{r\ge0}t^{-r}K_{n,U}t^{r}$. The last condition means
that $\Ad(t)$ strictly contracts $\Lie U$; hence $\Ad(t)^{-1}$
strictly contracts $\Lie\bar{U}$ and $\bar{U}=\bigcup_{r\ge0}t^{r}K_{n,\bar{U}}t^{-r}$,
so $t$ is strongly positive in the sense of \cite[\S3.2]{MR1810135}.
Put $e=e_{K_{n}}$, $\mathcal{A}=e\Hk(G)e$ and $h=e*\delta_{t}*e\in\mathcal{A}$.
The algebra $\mathcal{A}$ has the unit $e$. For a representation
$V$ we write $T$ for the action of $h$ on $V^{K_{n}}$ and $j_{U}:V\to V_{U}$
for the Jacquet map.
\begin{lemma}[The nilpotent part]
\label{lem:jacquet-kernel} For every representation $V$ and $r\ge0$
we have $T^{r}v=t^{r}e_{t^{-r}K_{n,U}t^{r}}v$ for $v\in V^{K_{n}}$.
Consequently $\ker(j_{U}|_{V^{K_{n}}})=\bigcup_{r\ge0}\ker T^{r}$. 
\end{lemma}

\begin{proof}
Apply \eqref{eq:positive-averaging} with $J=K_{n}$ and iterate,
using that the groups $t^{-r}K_{n,U}t^{r}$ increase with $r$.
They exhaust $U$, and a vector has zero Jacquet image precisely when
some compact open subgroup of $U$ averages it to zero.
\end{proof}

\begin{lemma}[Qualitative stabilization]
\label{lem:qualitative} There exist an idempotent $p\in\mathcal{A}$
and an element $d\in p\mathcal{A}p$, both commuting with $h$, and
an integer $c_{0}\ge1$ such that $hd=p$ and $h^{c_{0}}(e-p)=0$.
Consequently, for every representation $V$, $V^{K_{n}}=pV^{K_{n}}\oplus(e-p)V^{K_{n}}$,
$T$ is invertible on $pV^{K_{n}}$ and nilpotent on $(e-p)V^{K_{n}}$,
and $(e-p)V^{K_{n}}=\ker(j_{U}|_{V^{K_{n}}})$. 
\end{lemma}

\begin{proof}
Apply \cite[Theorem~1]{MR1810135} to $X=\C[G/K_{n}]=\Hk(G)e$. Its
positive Hecke operator is a nonzero scalar multiple of $h$, so its
powers have the same kernels and images. On $X^{K_{n}}=\mathcal{A}$
this gives, for some $c_{0}\ge1$, 
\[
\mathcal{A}=\ker L^{c_{0}}_{h}\oplus\im L^{c_{0}}_{h},
\]
with $L_{h}$ invertible on the image. Both summands are right ideals.
The projection onto the image and the inverse of $L_{h}$ there, extended
by zero on the kernel, are right $\mathcal{A}$-linear endomorphisms
of the unital right module $\mathcal{A}$. They are consequently left
multiplications by elements $p,d\in\mathcal{A}$. The stable decomposition
gives exactly the asserted relations. These relations hold on every
left $\mathcal{A}$-module; the last assertion follows from \cref{lem:jacquet-kernel}. 
\end{proof}

\begin{proposition}[Reduction to induced representations]
\label{prop:reduction} Let $N\ge1$ and suppose that for every parabolic
subgroup $Q=LN_{Q}$ of $G$ and every irreducible supercuspidal representation
$\sigma$ of $L$ (including $Q=G$) we have 
\[
h^{N}\ker\bigl((i^{G}_{Q}\sigma)^{K_{n}}\to(i^{G}_{Q}\sigma)_{U}\bigr)=0.
\]
Then $h^{N}(e-p)=0$. Consequently, for every representation $V$
and every $r\ge N$, $\ker T^{r}=\ker T^{N}=(e-p)V^{K_{n}}$ and $\im T^{r}=\im T^{N}=pV^{K_{n}}$. 
\end{proposition}

\begin{proof}
Let $a=h^{N}(e-p)\in\mathcal{A}$. By assumption and \cref{lem:qualitative},
$a$ annihilates the $K_{n}$-invariants of every $i^{G}_{Q}\sigma$,
and since $a=eae$, it annihilates $i^{G}_{Q}\sigma$. Every irreducible
representation of $G$ is a subquotient of such a representation \cite{BernsteinHarvard},
so that $a$ annihilates every irreducible representation. The separation
property of \cite[\S4.5]{MR2234120} now gives $a=0$. The rest follows
from the relations of \cref{lem:qualitative}. 
\end{proof}

\section{Linear stabilization}

\label{sec:stab} We retain the \levsupport{} \eqref{eq:levi-support}
and the adapted tower of \cref{subsec:adapted-filtrations}. For a
proper standard parabolic $P=MU$, let $A_{M}$ and $\epsilon_{P}$
be as in \cref{sec:intro}. The weights in $\Phi(A_{M},U)$ are the
restrictions of the relative roots in $\Phi_{U}$, so 
\[
\epsilon_{P}(a)=\min_{\alpha\in\Phi_{U}}\val(\alpha(a))\qquad(a\in A_{M}).
\]
Put $T_{n,a}=\pi(e_{K_{n}}*\delta_{a}*e_{K_{n}})$ on $V^{K_{n}}$.
\begin{theorem}[Stabilization for sufficiently contracting elements]
\label{thm:stab} There exists $c>0$, depending only on the fixed
field, group, filtration, and Levi support bounds, such that for every
proper standard $P=MU$, every $n\ge1$, every smooth representation
$V$, and every $a\in A_{M}$ satisfying $\epsilon_{P}(a)\ge cn$,
\[
V^{K_{n}}=\ker T_{n,a}\oplus\im T_{n,a}.
\]
The Jacquet map vanishes on the kernel and identifies the image with
$(V_{U})^{K_{n,M}}$. The operator is invertible on its image. For
fixed $P$ and $n$, these summands are independent of $a$ in the
indicated range. 
\end{theorem}

\begin{proof}[Proof of \cref{thm:stab}]
Let $P=MU$ be standard and $a\in A_{M}$ with $\epsilon_{P}(a)\ge cn$,
where $c$ will be chosen below. By \cref{subsec:adapted-filtrations},
the group $K_{n}$ satisfies the hypotheses of \cref{sec:smooth}.
Every parabolic is conjugate to a standard one, and conjugation does
not change the isomorphism class of the resulting induced representation.
By \cref{prop:reduction}, it suffices to kill the Jacquet kernel
in $(i^{G}_{Q}\sigma)^{K_{n}}$ for standard $Q=LN_{Q}$ and supercuspidal
$\sigma$.

Put $D'_{P}=\max_{\alpha\in\Phi_{U}}\val(\alpha(t_{P}))$. By
\eqref{eq:root-depths} and \eqref{eq:root-translation}, there is a constant
$c_{P}\ge0$ such that $b(K\cap U)b^{-1}\subseteq K_{n,U}$ whenever
$\epsilon_{P}(b)\ge n+c_{P}$. Indeed, conjugation increases each root depth
by at least $\epsilon_{P}(b)$, and the target depths are $n-1-\alpha(x)$,
including those for doubled roots. Taking $b=at_{P}^{-m}$ gives
\begin{equation}
a^{-1}K_{n,U}a\supset U[m]\quad\text{if}\quad
\epsilon_{P}(a)\ge n+c_{P}+D'_{P}m,\label{eq:final-containment}
\end{equation}
since $\epsilon_{P}(at_{P}^{-m})\ge\epsilon_{P}(a)-mD'_{P}$.

Choose an integer $C$ as in \cref{prop:induced} and put $m=Cn$.
For $f\in(i^{G}_{Q}\sigma)^{K_{n}}$ with zero Jacquet image, \cref{lem:jacquet-kernel}
and \eqref{eq:final-containment} give 
\[
T_{n,a}f=ae_{a^{-1}K_{n,U}a}e_{U[m]}f=0.
\]
Taking $c\ge 1+c_{P}+D'_{P}C$, and maximizing over the finitely
many standard parabolics, proves the stable decomposition and kernel
assertion by \cref{prop:reduction}, with exponent $1$. The same
reduction shows that $T_{n,a}$ is invertible on its image. The identification
with Jacquet invariants and independence of $a$ are the standard
properties of this stable summand \cite[Theorem~1 and \S3.5]{MR1810135}. 
\end{proof}

\begin{corollary}[Stabilization of powers]
\label{cor:stab-powers} With the hypotheses and constant of \cref{thm:stab},
let $t\in A_{M}$ satisfy $\epsilon_{P}(t)>0$. For $n\ge1$, put
$N_{n}(t)=\lceil cn/\epsilon_{P}(t)\rceil$. For every integer $r\ge N_{n}(t)$,
\[
V^{K_{n}}=\ker T^{\,r}_{n,t}\oplus\im T^{\,r}_{n,t}.
\]
The Jacquet map vanishes on the kernel and identifies the image with
$(V_{U})^{K_{n,M}}$. Both summands are independent of $r$ in this
range, and $T_{n,t}$ is invertible on the image. 
\end{corollary}

\begin{proof}
By \cref{lem:jacquet-kernel}, $T^{\,r}_{n,t}=T_{n,t^{r}}$. Since
$\epsilon_{P}(t^{r})=r\epsilon_{P}(t)$, \cref{thm:stab} applied
to $t^{r}$ gives the asserted summands and Jacquet isomorphism.
On the common image, $T_{n,t}$ is invertible because its $r$-th power is.
\end{proof}

\begin{remark*}[Sharpness of the stabilization bound]
 Fix a proper standard parabolic $P=MU$, a contracting $t\in A_{M}$
as in \cref{cor:stab-powers}, and $\pi\in\Irrsc G$. Choose a nonzero
$K_{b}$-fixed vector $u\in\pi$ and an integer $D\ge1$ as in
\eqref{eq:conj-shift} for $t$. For $n\ge b$ put $m=\lfloor(n-b)/D\rfloor$.
Since $t^{\pm m}K_{n}t^{\mp m}\subseteq K_{b}$, both $\pi(t^{-m})u$
and $\pi(t^{m})u$ are $K_{n}$-fixed. Hence
\[
T_{n,t}^{2m}\pi(t^{-m})u=T_{n,t^{2m}}\pi(t^{-m})u=\pi(t^{m})u\ne0.
\]
Since $\pi_{U}=0$, \cref{cor:stab-powers} implies that $T_{n,t}$ is
nilpotent on $\pi^{K_{n}}$. Its stabilization exponent is therefore
at least $2m+1$, which grows linearly in $n$ even for fixed $\pi$ and $t$.
\end{remark*}
\begin{remark}[Dependence on parameters]
\label{rem:parameters} Consider a family of bounded absolute rank
with a parameter $\eta\ge1$. Suppose the Levi support radii, the
stratum radius in \cref{lem:stratum}, and the depth bounds in \cref{eq:conj-shift,eq:levi-comparison,eq:unipotent-shift,eq:final-containment}
have uniformly bounded slopes and additive terms $O(\eta)$. For an
embedded family of Levi subgroups in characteristic zero, \eqref{eq:family-matrix-support}
supplies the matrix support estimate. The comparisons \eqref{eq:matrix-adjoint-balls}
and \eqref{eq:levi-comparison} transfer it to the adapted tower.
Under the stated bounds the conclusion of \cref{thm:stab} holds whenever
$\epsilon_{P}(a)\ge C(n+\eta)$, with $C$ uniform in the family.
Indeed, retaining the additive terms in the proof of \cref{prop:induced}
gives $m_{0}=0$ and $m_{j}\le A(n+m_{j-1})+B\eta$ with uniform $A,B$.
The number of steps is at most $\abs W$, uniformly bounded in the family,
so the averaging radius is $O(n+\eta)$. The
claim follows from \eqref{eq:final-containment}. For a contracting
$t$, this gives stabilization of $T^{\,r}_{n,t}$ for $r\ge\lceil C(n+\eta)/\epsilon_{P}(t)\rceil$. 
\end{remark}

\begin{remark}[Use of qualitative stabilization]
\label{rem:qualitative-role} Qualitative stabilization supplies
the idempotent of \cref{lem:qualitative}. No estimate for its support
or for the qualitative exponent enters the bound. 
\end{remark}

\begin{remark}[Other parabolics]
\label{rem:all-parabolics} Every parabolic containing $\adtor$
is conjugate to a standard one by an element of $N_{G}(\adtor)\cap K_{0}$.
These representatives normalize every $K_{n}$. The same constant
therefore works for these parabolics, in particular for the opposite
parabolic $\bar{P}$ and the element $t^{-1}$, with $\epsilon_{\bar{P}}(t^{-1})=\epsilon_{P}(t)$. 
\end{remark}

\section{Canonical lifts and the asymptotic pairing}

\label{sec:pairing} Let $\bar{P}=M\bar{U}$ be opposite to $P=MU$,
put $\Delta_{P}(m)=|\det(\Ad(m)|_{\Lie U})|$, and write $r_{P}(V)=\Delta^{-1/2}_{P}\otimes V_{U}$.
Let $j_{U},j_{\bar{U}}$ be the Jacquet maps. 
\begin{theorem}[A linear threshold for asymptotics]
\label{thm:adj} Under the hypotheses of \cref{thm:stab}, there
is a unique bilinear pairing $B_{P}$ between $r_{P}(V)$ and $r_{\bar{P}}(V^{\vee})$
such that 
\[
\langle\pi(a)v,v^{\vee}\rangle=\Delta_{P}(a)^{1/2}B_{P}\bigl(r_{P}(\pi)(a)j_{U}(v),j_{\bar{U}}(v^{\vee})\bigr)
\]
for every $n\ge1$, $v\in V^{K_{n}}$, $v^{\vee}\in(V^{\vee})^{K_{n}}$
and $a\in A_{M}$ with $\epsilon_{P}(a)\ge cn$, where $c$ is the
constant of \cref{thm:stab}. 
\end{theorem}

We use the canonical splitting and duality of \cite[Theorem~1, Proposition~7, and \S5]{MR1810135}.
For admissible representations these constructions are due to Casselman
\cite{CassNotes}; the extension to all smooth representations is
due to Bernstein \cite{BernsteinHarvard}. Let $t$ and $N_{n}(t)$
be as in \cref{cor:stab-powers}. Put 
\[
T=T_{n,t},\qquad V_{n,0}=\ker T^{N_{n}(t)},\qquad V_{n,+}=\im T^{N_{n}(t)}.
\]
For the opposite parabolic and $V^{\vee}$, use the transpose operator
$\bar{T}=\pi^{\vee}(e_{K_{n}}*\delta_{t^{-1}}*e_{K_{n}})$ and denote
its stable summands by $V^{\vee}_{n,0}$ and $V^{\vee}_{n,+}$. Here
\begin{equation}
(V^{\vee})^{J}=(V^{J})^{*}\quad\text{for every compact open }J,\label{eq:fixed-duality}
\end{equation}
with the full algebraic dual on the right: a functional on $V^{J}$
extends by composition with $e_{J}$. No admissibility is required.
\begin{lemma}[Canonical lifts]
\label{lem:lift} The Jacquet maps restrict to isomorphisms 
\begin{equation}
V_{n,+}\xrightarrow{\sim}(V_{U})^{K_{n,M}},\qquad V^{\vee}_{n,+}\xrightarrow{\sim}((V^{\vee})_{\bar{U}})^{K_{n,M}}.\label{eq:jacquet-iso}
\end{equation}
Writing $t_{U}$ for the unnormalized action on $V_{U}$, the inverse
$s_{P,n}$ of the first map is 
\[
s_{P,n}(\xi)=T^{N}v,\qquad N\ge N_{n}(t),\quad v\in V^{K_{n}},\quad j_{U}(v)=t^{-N}_{U}\xi.
\]
For $n'\ge n$ and $\xi\in(V_{U})^{K_{n,M}}$, $e_{K_{n}}s_{P,n'}(\xi)=s_{P,n}(\xi)$. 
\end{lemma}

\begin{proof}
The isomorphisms and level compatibility are \cite[Theorem~1 and Proposition~7]{MR1810135}.
By \cref{cor:stab-powers}, their stable summands are already obtained
at exponent $N_{n}(t)$. The displayed formula follows from $j_{U}T=t_{U}j_{U}$
and $T^{N}V_{n,0}=0$. 
\end{proof}

Let $s^{\vee}_{\bar{P},n}$ denote the inverse of the second map in
\eqref{eq:jacquet-iso}.
\begin{proposition}[The canonical pairing]
\label{prop:pairing} There is a natural $M$-invariant pairing $B_{P}$
identifying $r_{\bar{P}}(V^{\vee})$ with the smooth contragredient
of $r_{P}(V)$. On $K_{n,M}$-invariants it is given by 
\[
B_{P}(\xi,\eta)=\langle s_{P,n}(\xi),s^{\vee}_{\bar{P},n}(\eta)\rangle.
\]
The stable summands for $P$ and $\bar{P}$ are in duality and the
mixed pairings vanish. 
\end{proposition}

\begin{proof}
This is \cite[Theorem~2 and Lemmas~6--7]{MR1810135}. The normalization
by $\Delta^{-1/2}_{P}$ and $\Delta^{-1/2}_{\bar{P}}=\Delta^{1/2}_{P}$
cancels in the pairing. 
\end{proof}

\begin{proof}[Proof of \cref{thm:adj}]
If $\epsilon_{P}(a)\ge cn$, \cref{thm:stab} identifies $\im T_{n,a}$
with $V_{n,+}$. The mixed pairings vanish, so 
\begin{equation}
\langle\pi(a)v,v^{\vee}\rangle=\langle T_{n,a}v,v^{\vee}\rangle=B_{P}\bigl(a_{U}j_{U}(v),j_{\bar{U}}(v^{\vee})\bigr).\label{eq:unnormalized-pairing}
\end{equation}
Since $a_{U}=\Delta_{P}(a)^{1/2}r_{P}(\pi)(a)$, this gives the asserted
identity. It determines $B_{P}$ on each pair of $K_{n,M}$-fixed spaces,
since for $\epsilon_{P}(a)\ge cn$, the action of $r_{P}(\pi)(a)$ on the first
space is invertible and both Jacquet maps in \cref{lem:lift} are surjective.
Every pair of vectors is fixed at some common level.
\end{proof}

\section{Second adjointness at finite level}

\label{sec:adjunction}

Bernstein's second adjointness \cite[Theorem~3]{MR1810135} gives
a natural isomorphism 
\[
\Hom_{G}(i_{P}W,V)\simeq\Hom_{M}(W,r_{\bar{P}}V)
\]
for smooth representations $W$ of $M$ and $V$ of $G$. A finite-level
manifestation of this isomorphism is given by \cref{prop:adjunction}
below. Use the usual model of $i_{P}W$, with $f(mug)=\Delta_{P}(m)^{1/2}\sigma(m)f(g)$
and right translation by $G$, where $\sigma$ denotes the action
on $W$. Fix a Haar measure $d\bar{u}$ on $\bar{U}$. For $w\in W^{K_{n,M}}$,
let $f_{n,w}$ be the right-$K_{n}$-invariant section supported on
$PK_{n}$ with 
\[
f_{n,w}(1)=\operatorname{vol}_{\bar{U}}(K_{n,\bar{U}})^{-1}w.
\]
We normalize second adjointness by $\alpha(w)=j_{\bar{U}}(\Phi(f_{n,w}))$,
as in \cref{prop:adjunction}. Replacing $d\bar{u}$ by $c\,d\bar{u}$,
$c>0$, multiplies the map $\Phi\mapsto\alpha$ by $c^{-1}$. The
normalized averaging operators, canonical lifts, and quantitative
bounds are independent of this choice.
\begin{proposition}[Corresponding morphisms]
\label{prop:adjunction} Let $\Phi:i_{P}W\to V$ and $\alpha:W\to r_{\bar{P}}V$
correspond under second adjointness. For $w\in W^{K_{n,M}}$, 
\[
\Phi(f_{n,w})=s_{\bar{P},n}(\alpha(w)),\qquad\alpha(w)=j_{\bar{U}}\bigl(\Phi(f_{n,w})\bigr),
\]
where $s_{\bar{P},n}$ is the canonical lift of \cref{lem:lift} for
$\bar{P}$. These formulas are compatible under averaging between
levels. For $t\in A_{M}$ contracting $U$, take $N_{n}(t)$ as in
\cref{cor:stab-powers}. On $V^{K_{n}}$ the operator for the opposite
parabolic is 
\[
T_{n,t^{-1}}=\pi(e_{K_{n}}*\delta_{t^{-1}}*e_{K_{n}}).
\]
Since $\epsilon_{\bar{P}}(t^{-1})=\epsilon_{P}(t)$, the Jacquet map
restricts to an isomorphism 
\[
j_{\bar{U}}:\im T^{N_{n}(t)}_{n,t^{-1}}\xrightarrow{\;\sim\;}(V_{\bar{U}})^{K_{n,M}}
\]
with inverse $s_{\bar{P},n}$. In particular, $\Phi(f_{n,w})$ is
the unique vector in $\im T^{N_{n}(t)}_{n,t^{-1}}$ whose Jacquet
image is $\alpha(w)$. 
\end{proposition}

\begin{proof}
Bushnell's definition of second adjointness \cite[\S6.1, Lemma~8 and (6.1)]{MR1810135}
gives $\alpha(w)=j_{\bar{U}}(\Phi(f_{n,w}))$ with our normalization.
Put $a=t^{-1}$ and $C=K_{n,\bar{U}}$. Since $PK_{n}=PC$, the translate
$af_{n,w}$ is supported on $P(aCa^{-1})$. Its restriction to $\bar{U}$ is
\[
\operatorname{vol}(C)^{-1}\Delta_{P}(a)^{1/2}\mathbf{1}_{aCa^{-1}}\,\sigma(a)w.
\]
By \eqref{eq:positive-averaging} for $\bar{P}$, $T_{n,a}f_{n,w}=e_{C}af_{n,w}$.
As $aCa^{-1}\subseteq C$ has volume $\Delta_{P}(a)^{-1}\operatorname{vol}(C)$,
averaging gives $T_{n,a}f_{n,w}=\Delta_{P}(a)^{-1/2}f_{n,\sigma(a)w}$.
Since $a$ is central in $M$, this action on the standard sections is invertible, and
$G$-equivariance puts $\Phi(f_{n,w})$ in $\im T_{n,a}^{N_{n}(t)}$.
By \cref{lem:lift} for $\bar{P}$, it is the unique lift there of $\alpha(w)$,
and these lifts are compatible under averaging.
\end{proof}

\raggedbottom  \bibliographystyle{amsalpha}
\bibliography{../../Bibfiles/all}
\end{document}